# On $n$-generalized commutators and Lie ideals of rings

Peter V. Danchev[♮] and Tsiu-Kwen Lee[‡]

Institute of Mathematics and Informatics, Bulgarian Academy of Sciences[♮]
"Acad. G. Bonchev" str., bl. 8, 1113 Sofia, Bulgaria
emails: danchev@math.bas.bg; pvdanchev@yahoo.com
and
Department of Mathematics, National Taiwan University[‡]
Taipei 106, Taiwan
email: tklee@math.ntu.edu.tw

**Abstract**

Let $R$ be an associative ring. Given a positive integer $n \geq 2$, for $a_1, \ldots, a_n \in R$ we define $[a_1, \ldots, a_n]_n := a_1a_2 \cdots a_n - a_na_{n-1} \cdots a_1$, the $n$-generalized commutator of $a_1, \ldots, a_n$. By an $n$-generalized Lie ideal of $R$ (at the $(r+1)$-th position with $r \geq 0$) we mean an additive subgroup $A$ of $R$ satisfying $[x_1, \ldots, x_r, a, y_1, \ldots, y_s]_n \in A$ for all $x_i, y_j \in R$ and all $a \in A$, where $r + s = n - 1$. In the paper we study $n$-generalized commutators of rings and prove that if $R$ is a noncommutative prime ring and $n \geq 3$, then every nonzero $n$-generalized Lie ideal of $R$ contains a nonzero ideal. Therefore, if $R$ is a noncommutative simple ring, then $R = [R, \ldots, R]_n$. This extends a classical result due to Herstein (Portugal. Math., 1954). Some generalizations and related questions on $n$-generalized commutators and their relationship with noncommutative polynomials are also discussed.

Corresponding author: Tsiu-Kwen Lee

# 1 Introduction

Throughout, $R$ always denotes an associative ring, not necessarily with unity, with center $Z(R)$. Given $x, y \in R$, let $[x, y] := xy - yx$, the additive commutator (or the Lie product) of $a$ and $b$. For additive subgroups $A, B$ of $R$, let $AB$ (resp. $[A, B]$) denote the additive subgroup of $R$ generated by all elements $ab$ (resp. $[a, b]$) for $a \in A$ and $b \in B$. An additive subgroup $L$ of $R$ is called a Lie ideal of $R$ if $[L, R] \subseteq L$. Given two Lie ideals $L_1$ and $L_2$ of $R$, it is clear that both $L_1L_2$ and $[L_1, L_2]$ are Lie ideals of $R$. A Lie ideal $L$ of $R$ is called proper if $[M, R] \subseteq L$ for some nonzero ideal $M$ of $R$. Given a nonempty subset $T$ of $R$, let $\mathcal{I}(T)$ denote the ideal of $R$ generated by $T$. A ring $R$ is called prime if, for $a, b \in R$, $aRb = 0$ implies that either $a = 0$ or $b = 0$.

Referring to [26], by a commutator ring we mean a ring $R$ satisfying $R = [R, R]$. In 1956 Kaplansky proposed twelve problems in ring theory (see [16]). One of these was whether there is a division ring $D$ satisfying $D = [D, D]$. Harris constructed the first example of division commutator rings (see [13]). Related researches are still quite active up to now (see [26, 3, 30] *et al.*).

For $a, b, c \in R$, we let $[a, b, c] := abc - cba$, which is called the generalized commutator of $a, b$ and $c$. Given additive subgroups $A_1, A_2, A_3$ of $R$, let $[A_1, A_2, A_3]$ denote the additive subgroup of $R$ generated by all elements $[a_1, a_2, a_3]$ for $a_i \in A_i$, $i = 1, 2, 3$. In 1954 Herstein initiated the study of generalized commutators (see [14]). An independent work dealing with generalized commutators of matrices over fields was due to Taussky (see [29]). Herstein proved that $[R, R, R]$ is an ideal of $R$ (see [14, Theorem 3]) and is equal to the whole ring $R$ if $R$ is a noncommutative simple ring (see [14, Theorem 4]). Further, if $R$ is a simple Artinian ring, then every element of $R$ is a sum of three generalized commutators (see [14, Theorem 6]). Generalized commutators also naturally appear in analysis. A classical analytic result of Brown and Pearcy [4], which states that, in a $\mathbb{C}$-algebra of bounded operators on a Hilbert space of dimension $\geq 2$, every operator is a generalized commutator (see also [29]). Relevant researches have been brought up again in recent years. See [11, 17] and the references therein. In the paper we study generalized commutators in a more generalized form.

**Definition 1.** Given a positive integer $n \geq 2$, for $a_1, \ldots, a_n \in R$ we define

$$[a_1, \ldots, a_n]_n := a_1a_2 \cdots a_n - a_na_{n-1} \cdots a_1,$$

which is called the $n$-generalized commutator of $a_1, \ldots, a_n$.

Given additive subgroups $A_1, A_2, \ldots, A_n$ of $R$, let $[A_1, A_2, \ldots, A_n]_n$ denote the additive subgroup of $R$ generated by all elements $[a_1, a_2, \ldots, a_n]_n$ for $a_i \in A_i$, $i = 1, \ldots, n$. If $n = 2$, then $[a_1, a_2]_2 = [a_1, a_2] = a_1a_2 - a_2a_1$, the ordinary additive commutator of $a_1$ and $a_2$. If $n = 3$, then $[a_1, a_2, a_3]_3$ coincides with $[a_1, a_2, a_3]$, the generalized commutator of $a_1, a_2, a_3$. The purpose of this paper is to study the following:

- Properties of $n$-generalized commutators (see §2 and §3).
- $n$-Generalized commutator rings (see §4).
- $n$-Generalized Lie ideals as a generalization of Lie ideals (see §5, §6, and §7).
- Some generalizations connected to noncommutative polynomials (see §8).

# 2 $n$-Generalized commutators

We begin with an observation. Given an odd positive integer $n > 1$, the following two equalities hold in an arbitrary ring $R$:

$$\begin{aligned} &[a_1, a_2, \ldots, a_n]_n r \\ &= [a_1, a_2, \ldots, a_n r]_n - [a_1, a_2, \ldots, r a_{n-1}, a_n]_n \\ &\qquad + [a_1, a_2, \ldots, a_{n-2} r, a_{n-1}, a_n]_n - \cdots + [a_1 r, a_2, \ldots, a_n]_n \end{aligned} \tag{1}$$

and

$$\begin{aligned} &r[a_1, a_2, \ldots, a_n]_n \\ &= [a_1, a_2, \ldots, r a_n]_n - [a_1, a_2, \ldots, a_{n-1} r, a_n]_n \\ &\qquad + [a_1, a_2, \ldots, r a_{n-2}, a_{n-1}, a_n]_n - \cdots + [r a_1, a_2, \ldots, a_n]_n. \end{aligned} \tag{2}$$

Therefore, the following gives a generalization of [14, Theorem 3].

**Theorem 2.1.** *Let $R$ be a ring. Then $[R, \ldots, R]_{2n+1}$ is an ideal of $R$ for $n \geq 1$.*

Clearly, Theorem 2.1 has a more generalized form: If $I_1, \ldots, I_{2n+1}$ are ideals of $R$, then $[I_1, \ldots, I_{2n+1}]_{2n+1}$ is an ideal of $R$ for $n \geq 1$. Theorem 2.1 serves as the starting point for understanding this nature of $[R, \ldots, R]_n$. A natural question is to ask whether $[R, \ldots, R]_n$ is an ideal of $R$ for even $n$. Clearly, it is not in general true for $n = 2$. We begin with some basic observations.

**Lemma 2.2.** *Let $R$ be a ring with a Lie ideal $L$. Then the following hold:*
*(i) If $I$ is an ideal of $R$, then $IL$ and $LI$ are ideals of $R$.*
*(ii) $L + LR = \mathcal{I}(L) = L + RL$.*

*Proof.* (i) Clearly, $RIL \subseteq IL$. Also, $ILR \subseteq I([L, R] + RL) \subseteq IL + IRL \subseteq IL$. Therefore, $RL$ is an ideal of $R$. Similarly, $LI$ is an ideal of $R$.

(ii) Note that $RL \subseteq [L, R] + LR \subseteq L + LR$. By (i), $RLR \subseteq LR$. Therefore, $\mathcal{I}(L) = L + LR + RL + RLR \subseteq L + LR$ and so $\mathcal{I}(L) = L + LR$. The another case has the same argument. □

**Proposition 2.3.** *Let $R$ be a ring. Then the following hold:*

*(i)* $[R,\ldots,R]_n \subseteq \mathcal{I}([R,R])$ *for* $n \geq 2$.

*(ii)* $[R,\ldots,R]_{2n-1} \subseteq [R,R,R]$ *for* $n \geq 2$.

*(iii)* $[R,\ldots,R]_{2n} \subseteq [R,R] + [R,\ldots,R]_{2n-1} \subseteq \mathcal{I}([R,R])$ *for* $n \geq 2$.

*(iv)* $[R,R] + [R,R,R] = \mathcal{I}([R,R])$.

*(v)* $[R,\ldots,R]_{2n} \subseteq [R,R,R,R] + [R,\ldots,R]_{2n-3} \subseteq [R,R,R] + [R,R,R,R]$ *for* $n \geq 3$.

*(vi) If* $R = R^2$, *then* $[R,\ldots,R]_k + [R,\ldots,R]_{k+1} = \mathcal{I}([R,R])$ *for* $k \geq 2$.

*(vii) If* $R = R^2$, *then* $[R,R] + [R,\ldots,R]_{2n-1} = \mathcal{I}([R,R])$ *for* $n \geq 2$.

*Proof.* (i) The case $n = 2$ is trivial. In view of [14, Lemma 1], $[R,R,R] \subseteq \mathcal{I}([R,R])$. Assume that $n > 3$. By induction on $n$, we assume that $[R,\ldots,R]_{n-1} \subseteq \mathcal{I}([R,R])$. Let $a_1,\ldots,a_n \in R$. Then

$$\begin{aligned}&[a_1,\ldots,a_n]_n\\ &= [a_1,\ldots,a_{n-1}]_{n-1}a_n - \big[a_n, a_{n-1}a_{n-2}\cdots a_1\big] \in [R,\ldots,R]_{n-1}a_n + [R,R]\\ &\subseteq \mathcal{I}([R,R]).\end{aligned}$$

Therefore, $[R,\ldots,R]_n \subseteq \mathcal{I}([R,R])$.

(ii) It is trivial for $n = 2$. Assume that $n \geq 3$. Applying the inductive argument, we assume that $[R,\ldots,R]_{2n-3} \subseteq [R,R,R]$. Let $a_1,\ldots,a_{2n-1} \in R$. Then

$$\begin{aligned}&[a_1,\ldots,a_{2n-1}]_{2n-1}\\ &= [a_1a_2a_3, a_4,\ldots,a_{2n-1}]_{2n-3} + a_{2n-1}a_{2n-2}\cdots a_4[a_1,a_2,a_3]\\ &\in [R,\ldots,R]_{2n-3} + [R,R,R],\end{aligned}$$

since, by Theorem 2.1, $[R,R,R]$ is an ideal of $R$. Therefore,

$$[R,\ldots,R]_{2n-1} \subseteq [R,\ldots,R]_{2n-3} + [R,R,R] = [R,R,R].$$

(iii) Applying the computation in (i), we get

$$[R,\ldots,R]_{2n} \subseteq [R,R] + [R,\ldots,R]_{2n-1}R \subseteq [R,R] + [R,\ldots,R]_{2n-1}$$

since, by Theorem 2.1, $[R,\ldots,R]_{2n-1}$ is an ideal of $R$. By (i), $[R,\ldots,R]_{2n-1} \subseteq \mathcal{I}([R,R])$. Hence $[R,\ldots,R]_{2n} \subseteq [R,R] + [R,\ldots,R]_{2n-1} \subseteq \mathcal{I}([R,R])$, as desired.

(iv) By (i), it is clear that $[R,R] + [R,R,R] \subseteq \mathcal{I}([R,R])$. Let $a,b,x \in R$. Then

$$[a,b]x = [a,b,x] + [x,ba].$$

This means that $[R,R]R \subseteq [R,R,R] + [R,R]$. In view of (ii) of Lemma 2.2,

$$\mathcal{I}([R,R]) = [R,R] + [R,R]R \subseteq [R,R,R] + [R,R].$$

Hence $\mathcal{I}([R,R]) = [R,R,R] + [R,R]$.

(v) Let $a_1, \ldots, a_{2n} \in R$. Since $[R, \ldots, R]_{2n-3}$ is an ideal of $R$ and is contained in $[R,R,R]$ by (ii), we have

$$\begin{aligned} & [a_1, a_2, \ldots, a_{2n}]_{2n} \\ = & [a_1 a_2 \cdots a_{2n-3}, a_{2n-2}, a_{2n-1}, a_{2n}]_4 + a_{2n} a_{2n-1} a_{2n-2} [a_1, a_2, \ldots, a_{2n-3}]_{2n-3} \\ \in & [R,R,R,R] + [R, \ldots, R]_{2n-3} \subseteq [R,R,R,R] + [R,R,R]. \end{aligned}$$

Hence $[R, \ldots, R]_{2n} \subseteq [R,R,R,R] + [R, \ldots, R]_{2n-3} \subseteq [R,R,R] + [R,R,R,R]$.

(vi) **Case 1**: $k = 2n-1$, where $n \geq 2$. Let $a_1, \ldots, a_{2n} \in R$. Since $[R, \ldots, R]_{2n-1}$ is an ideal of $R$, we have

$$\begin{aligned} & [a_1 a_2 \cdots a_{2n-1}, a_{2n}] \\ = & [a_1, a_2, \ldots, a_{2n}]_{2n} - a_{2n} [a_1, a_2, \ldots, a_{2n-1}]_{2n-1} \\ \in & [R, \ldots, R]_{2n} + a_{2n} [R, \ldots, R]_{2n-1} \subseteq [R, \ldots, R]_{2n} + [R, \ldots, R]_{2n-1}. \end{aligned}$$

Therefore, $[R^{2n-1}, R] \subseteq [R, \ldots, R]_{2n} + [R, \ldots, R]_{2n-1}$. By $R = R^2$, we have $R = R^{2n-1}$. Hence $[R,R] = [R^{2n-1}, R] \subseteq [R, \ldots, R]_{2n} + [R, \ldots, R]_{2n-1}$.

We also have

$$\begin{aligned} & [a_{2n}, a_{2n-1}] a_{2n-2} \cdots a_2 a_1 \\ = & [a_{2n}, a_{2n-1}, \ldots, a_1]_{2n} + [a_1, \ldots, a_{2n-2}, a_{2n-1} a_{2n}]_{2n-1} \\ \in & [R, \ldots, R]_{2n} + [R, \ldots, R]_{2n-1}. \end{aligned}$$

Therefore, we have $[R,R]R = [R,R]R^{2n-2} \subseteq [R, \ldots, R]_{2n} + [R, \ldots, R]_{2n-1}$. In view of Lemma 2.2 and (i), we get

$$\mathcal{I}([R,R]) = [R,R] + [R,R]R \subseteq [R, \ldots, R]_{2n} + [R, \ldots, R]_{2n-1} \subseteq \mathcal{I}([R,R]).$$

Hence $[R, \ldots, R]_{2n} + [R, \ldots, R]_{2n-1} = \mathcal{I}([R,R])$.

**Case 2**: $k = 2n$, where $n \geq 1$. Let $a_1, \ldots, a_{2n+1} \in R$. Then

$$\begin{aligned} & [a_1, a_2] a_3 \cdots a_{2n} a_{2n+1} \\ = & [a_1, \ldots, a_{2n+1}]_{2n+1} - [a_2 a_1, a_3, \ldots, a_{2n}, a_{2n+1}]_{2n} \\ \in & [R, \ldots, R]_{2n+1} + [R, \ldots, R]_{2n}. \end{aligned}$$

Note that $R = R^2 = R^{2n-1}$. Therefore, $[R,R]R \subseteq [R, \ldots, R]_{2n+1} + [R, \ldots, R]_{2n}$. Similarly, $R[R,R] \subseteq [R, \ldots, R]_{2n+1} + [R, \ldots, R]_{2n}$. We then have

$$[R,R] = [R,R^2] \subseteq R[R,R] + [R,R]R \subseteq [R, \ldots, R]_{2n+1} + [R, \ldots, R]_{2n}.$$

By (ii) of Lemma 2.2 and (i), we have $[R,\ldots,R]_{2n+1}+[R,\ldots,R]_{2n}=\mathcal{I}([R,R])$.

(vii) In view of (iii) and (vi),

$$\begin{aligned}&\mathcal{I}([R,R])=[R,\ldots,R]_{2n}+[R,\ldots,R]_{2n-1}\\ &\subseteq [R,\ldots,R]_{2n}+[R,R]+[R,\ldots,R]_{2n-1}\\ &=[R,R]+[R,\ldots,R]_{2n-1}\subseteq \mathcal{I}([R,R]).\end{aligned}$$

This implies that $[R,R]+[R,\ldots,R]_{2n-1}=\mathcal{I}([R,R])$. □

Recall that, given a positive integer $n$, $R^n$ denotes the additive subgroup of $R$ generated by all elements $a_1a_2\cdots a_n$ for $a_i\in R$. Example (1) below shows that the inclusions in (i), (ii) and (iii) of Proposition 2.3 are proper. Example (2) below shows that the assumption that $R=R^2$ is essential to (vi) of Proposition 2.3.

Given a ring $R$ and $n$ a positive integer, let $\mathrm{M}_n(R)$ denote the $n\times n$ matrix ring over the ring $R$. If $1\in R$, we let $e_{ij}$'s, $1\le i,j\le n$, be the usual matrix units of $\mathrm{M}_n(R)$.

**Examples.** **(1)** Let $R:=\mathrm{M}_m(2\mathbb{Z})$, where $m\ge 2$. Then $[R,\ldots,R]_n\subsetneq\mathcal{I}([R,R])$, and $[R,\ldots,R]_{2n-1}\subsetneq[R,R,R]$ for any $n\ge 3$. Moreover, $[R,\ldots,R]_{2n}\subsetneq[R,R]+[R,\ldots,R]_{2n-1}$ for any $n\ge 2$.

In view of Proposition 2.3, it suffices to claim that $[R,R]\not\subseteq[R,\ldots,R]_n$ for $n\ge 3$. Then $4e_{12}=[2e_{12},2e_{22}]\in[R,R]$ but $[R,\ldots,R]_n\subseteq\mathrm{M}_m(8\mathbb{Z})$ since $n\ge 3$. This proves our claim. The other two cases have similar arguments. □

**(2)** Given a positive integer $m>1$, let $ST_m(K)$ denote the strict upper triangular subring of $\mathrm{M}_m(K)$, where $K$ is a unital ring. That is,

$$ST_m(K)=\{\sum_{1\le i<j\le m}a_{ij}e_{ij}\in\mathrm{M}_m(K)\mid a_{ij}\in K\ \forall i,j\}.$$

Set $R:=ST_m(K)$. Note that $R^m=0$ and, for $1\le k\le m-1$, we have

$$R^k=\{\sum_{1\le i;\ i+k\le j\le m}a_{ij}e_{ij}\in\mathrm{M}_m(K)\mid a_{ij}\in K\ \forall i,j\}.$$

A direct computation shows that $[R,\ldots,R]_k=R^k$ for $k\ge 2$. Moreover, we have

$$\{0\}=R^m\subsetneq R^{m-1}\subsetneq\cdots\subsetneq R^2\subsetneq R.$$

Therefore, for $m>4$, we have

$$[R,R,R]+[R,R,R,R]\subsetneq[R,R]=\mathcal{I}([R,R]).$$

□

Clearly, $R=R^2$ if $R$ is a unital ring. Continuing Proposition 2.3, we further study the nature of $[R,\ldots,R]_n$ for $n\ge 2$.

**Theorem 2.4.** *Let $R$ be a ring satisfying $R = R^2$. Then the following hold:*

*(i)* $[R, \ldots, R]_{2n-1} + [R, \ldots, R]_{2n-3} = [R, R, R]$ *for* $n \geq 3$.

*(ii) If* $[R, \ldots, R]_{2n}$ *is an ideal of* $R$*, then* $[R, \ldots, R]_{2n} = \mathcal{I}([R, R])$ *for* $n \geq 1$.

*Proof.* (i) We first claim that $[R, R, R] \subseteq [R, \ldots, R]_{2n-1} + [R, \ldots, R]_{2n-3}$ for $n \geq 3$. Let $a_1, \ldots, a_{2n-1} \in R$. In view of Theorem 2.1, $[R, \ldots, R]_{2n-3}$ is an ideal of $R$. Then

$$\begin{aligned}
&[a_1, a_2, a_3a_4 \cdots a_{2n-1}] \\
&= [a_1, a_2, \ldots, a_{2n-2}, a_{2n-1}]_{2n-1} + (a_{2n-1}a_{2n-2} \cdots a_3)a_2a_1 - (a_3a_4 \cdots a_{2n-1})a_2a_1 \\
&= [a_1, a_2, \ldots, a_{2n-2}, a_{2n-1}]_{2n-1} - [a_3, a_4, \ldots, a_{2n-2}, a_{2n-1}]_{2n-3}a_2a_1 \\
&\in [R, \ldots, R]_{2n-1} + [R, \ldots, R]_{2n-3}.
\end{aligned}$$

Since $R = R^{2n-3}$, we get $[R, R, R] \subseteq [R, \ldots, R]_{2n-1} + [R, \ldots, R]_{2n-3}$ for $n \geq 3$. On the other hand, by (ii) of Proposition 2.3, we have

$$[R, \ldots, R]_{2n-1} + [R, \ldots, R]_{2n-3} \subseteq [R, R, R]$$

for $n \geq 3$. Therefore, $[R, \ldots, R]_{2n-1} + [R, \ldots, R]_{2n-3} = [R, R, R]$.

(ii) It is trivial for $n = 1$. Suppose that $n > 1$.

**Step 1**: $R[R, R] \subseteq [R, \ldots, R]_{2n}$. Let $a_1, \ldots, a_{2n}, r \in R$. Since $[R, \ldots, R]_{2n}$ is an ideal of $R$, we have

$$\begin{aligned}
&[a_1, \ldots, a_{2n}]_{2n}r \\
&= [a_1, \ldots, a_{2n}r]_{2n} - [a_1, \ldots, ra_{2n-1}, a_{2n}]_{2n} + \cdots + [a_1, a_2r \ldots, a_{2n}]_{2n} \\
&\qquad + a_{2n}a_{2n-1} \cdots a_2[r, a_1] \in [R, \ldots, R]_{2n},
\end{aligned}$$

implying $a_{2n}a_{2n-1} \cdots a_2[r, a_1] \in [R, \ldots, R]_{2n}$. That is, $R^{2n-1}[R, R] \subseteq [R, \ldots, R]_{2n}$. Recall that $R = R^2$. We get $R^{2n-1} = R$ and so $R[R, R] \subseteq [R, \ldots, R]_{2n}$.

**Step 2**: $[R, \ldots, R]_{2n-1} \subseteq [R, \ldots, R]_{2n}$. Let $a_1, \ldots, a_{2n-2}, x, y \in R$. Then

$$\begin{aligned}
&[xy, a_1, \ldots, a_{2n-2}]_{2n-1} \\
&= [x, y, a_1, \ldots, a_{2n-2}]_{2n} - a_{2n-2} \cdots a_2a_1[x, y] \in [R, \ldots, R]_{2n},
\end{aligned}$$

since $a_{2n-2} \cdots a_2a_1[x, y] \in [R, \ldots, R]_{2n}$ by Step 1. Hence we have

$$[R^2, R, \ldots, R]_{2n-1} \subseteq [R, \ldots, R]_{2n}.$$

By $R = R^2$, we conclude that $[R, R, \ldots, R]_{2n-1} \subseteq [R, \ldots, R]_{2n}$, as desired.

It follows from both Step 2 and (vi) of Proposition 2.3 that

$$[R, \ldots, R]_{2n} = [R, R, \ldots, R]_{2n-1} + [R, \ldots, R]_{2n} = \mathcal{I}([R, R]),$$

as desired. □

The following is an immediate consequence of both (vi) of Proposition 2.3 and (ii) of Theorem 2.4.

**Corollary 2.5.** *Let $R$ be a ring satisfying $R = R^2$, $n \geq 2$ a positive integer. Then the following are equivalent:*

*(i)* $[R, \ldots, R]_{2n} = \mathcal{I}([R, R])$.

*(ii) Either* $[R, \ldots, R]_{2n-1} \subseteq [R, \ldots, R]_{2n}$ *or* $[R, \ldots, R]_{2n+1} \subseteq [R, \ldots, R]_{2n}$.

Similarly, we also have the following corollary.

**Corollary 2.6.** *Let $R$ be a ring satisfying $R = R^2$, $n \geq 2$ a positive integer. Then the following are equivalent:*

*(i)* $[R, \ldots, R]_{2n-1} = \mathcal{I}([R, R])$.

*(ii) Either* $[R, \ldots, R]_{2n} \subseteq [R, \ldots, R]_{2n-1}$ *or* $[R, \ldots, R]_{2n-2} \subseteq [R, \ldots, R]_{2n-1}$.

Up to now, the following problem keeps unknown.

**Problem 1.** Let $R$ be a ring satisfying $R = R^2$, and $n \geq 2$ a positive integer.

(i) Is $[R, \ldots, R]_{2n}$ an ideal of $R$?

(ii) Is $[R, \ldots, R]_{2n-1}$ equal to $\mathcal{I}([R, R])$?

We will answer Problem 1 affirmatively for unital rings, regular rings and rings generated by idempotents *et al.* in the next section.

# 3 Results on Problem 1

For (i) of Problem 1, we have the following observation.

**Proposition 3.1.** *Let $R$ be a ring satisfying $R = R^2$, $n$ a positive integer $n \geq 2$. Then the following are equivalent:*

*(i)* $[R, \ldots, R]_{2n}$ *is an ideal of* $R$.

*(ii)* $R[R, R] \subseteq [R, \ldots, R]_{2n}$.

*(iii)* $[R, R]R \subseteq [R, \ldots, R]_{2n}$.

*Proof.* "(i) $\Leftrightarrow$ (ii)": Set $L := [R, \ldots, R]_{2n}$. It follows from Theorem 2.4 that if $L$ is an ideal of $R$, then $L = \mathcal{I}([R, R])$ and hence $R[R, R] \subseteq L$. Conversely, assume that $R[R, R] \subseteq L$.

Let $a_1, \ldots, a_{2n}, r \in R$. By the proof of Step 1 in (ii) of Theorem 2.4,

$$[a_1, \ldots, a_{2n}]_{2n} r + a_{2n} a_{2n-1} \cdots a_2 [r, a_1] \in [R, \ldots, R]_{2n} = L.$$

Since $a_{2n}a_{2n-1}\cdots a_2[r,a_1] \in R[R,R] \subseteq L$, we have $[a_1,\ldots,a_{2n}]_{2n}r \in L$. That is, $LR \subseteq L$. Therefore, $L$ is a right ideal of $R$. Since $L$ is a Lie ideal of $R$, we have $RL \subseteq [R,L] + LR \subseteq L$. This proves that $L$ is an ideal of $R$.

By symmetry, we also have the same argument for the equivalence "(i) $\Leftrightarrow$ (iii)". □

**Lemma 3.2.** *Let $R$ be a ring. Suppose that $[R,\ldots,R]_k \subseteq [R,\ldots,R]_{k+1}$ for all $k \geq 2$. Then $[R,\ldots,R]_n = \mathcal{I}([R,R])$ for all $n \geq 3$.*

*Proof.* In view of Proposition 2.3, we have

$$[R,R] \subseteq [R,R,R] \subseteq [R,\ldots,R]_n \subseteq \mathcal{I}([R,R])$$

for $n \geq 3$. Since $[R,R,R]$ is an ideal of $R$, we get $\mathcal{I}([R,R]) \subseteq [R,R,R]$ and hence $[R,\ldots,R]_n = \mathcal{I}([R,R])$ for all $n \geq 3$. □

Let $R$ be an algebra over a field $F$. If $\operatorname{char} F \neq 2$, then $2R = R$. Before proving the next theorem, we need the following lemma.

**Lemma 3.3.** *Let $R$ be a ring. If $R = R^2$ and $2R = R$, then $R$ is equal to its additive subgroup generated by all square elements of $R$.*

*Proof.* We denote by $A$ the additive subgroup of $R$ generated by all square elements of $R$. Let $x,y,z \in A$. Then $xy + yx = (x+y)^2 - x^2 - y^2 \in A$. Therefore, we have

$$zxy - xyz = \big(y(zx) + (zx)y\big) - \big(x(yz) + (yz)x\big) \in A.$$

This implies that $2zxy = \big(z(xy) + (xy)z\big) + (zxy - xyz) \in A$. We thus get $2R^3 \subseteq A$. Since $2R = R$ and $R = R^2$, we conclude that $A = R$, as desired. □

**Theorem 3.4.** *Let $R$ be a ring satisfying both $R = R^2$ and $2R = R$. Then $[R,\ldots,R]_n = \mathcal{I}([R,R])$ for all $n \geq 3$.*

*Proof.* In view of Lemma 3.2, it suffices to claim that $[R,\ldots,R]_k \subseteq [R,\ldots,R]_{k+1}$ for all $k \geq 2$. Let $y, x_1, \ldots, x_{k-1} \in R$, where $k \geq 2$. We have

$$[y^2, x_1, \ldots, x_{k-1}]_k = [y, y, x_1, \ldots, x_{k-1}]_{k+1} \in [R,\ldots,R]_{k+1}.$$

That is,

$$[y^2, R, \ldots, R]_k \subseteq [R,\ldots,R]_{k+1}$$

for all $y \in R$. Since $R = R^2$ and $2R = R$, it follows from Lemma 3.3 that $R$ is equal to its additive subgroup generated by all square elements of $R$. We get $[R,R,\ldots,R]_k \subseteq [R,\ldots,R]_{k+1}$, as claimed. □

A ring $R$ is said to have the property (*) if given $x, y \in R$ there exists an element $c \in R$ such that $x, y \in cR$. Clearly, every unital ring has the property (*). Moreover, if $R$ has the property (*), then for finitely many $x_1, \ldots, x_n \in R$ there exists an element $c \in R$ such that $x_i \in cR$ for $i = 1, \ldots, n$. A ring $R$ is called regular if, given $a \in R$, there exists an element $a^- \in R$ satisfying $aa^-a = a$. In view of [10, Theorem 1.1] (it is true even when the ring $R$ has no unity), any regular ring $R$ has the property (*). Note that every ring $R$ having the property (*) always satisfies $R = R^2$.

**Proposition 3.5.** *Let $R$ be a ring having the property (*). Then $[R, R] + [R, R, R, R] = \mathcal{I}([R, R])$.*

*Proof.* In view of (i) of Proposition 2.3, $[R, R] + [R, R, R, R] \subseteq \mathcal{I}([R, R])$. Let $u, v, w \in R$. Since $R$ satisfies the property (*), there exists $a \in R$ such that $u, w \in aR$. Write $u = ax$ and $w = ay$ for some $x, y \in R$. Then

$$[u, v]w = [ax, v]ay = [a, x, va, y] + [y, vaxa] \in [R, R, R, R] + [R, R].$$

That is, $[R, R]R \subseteq [R, R, R, R] + [R, R]$. Therefore, it follows from (ii) of Lemma 2.2 that

$$\mathcal{I}([R, R]) = [R, R] + [R, R]R \subseteq [R, R, R, R] + [R, R],$$

implying that $\mathcal{I}([R, R]) = [R, R, R, R] + [R, R]$. □

The following corollary is an immediate consequence of both Proposition 3.5 and (ii) of Theorem 2.4.

**Corollary 3.6.** *Let $R$ be a ring having the property (*). Then $[R, R, R, R]$ is an ideal of $R$ if and only if $[R, R] \subseteq [R, R, R, R]$.*

We now answer Problem 1 for regular rings in the affirmative. A ring $R$ is said to have the property ($\sharp$) if, for any $x \in R$, there exists an idempotent $g \in R$ such that $x \in gRg$. Clearly, every ring having the property ($\sharp$) satisfies $R = R^2$.

**Theorem 3.7.** *Let $R$ be a ring satisfying the property ($\sharp$). Then $[R, \ldots, R]_n = \mathcal{I}([R, R])$ for all $n \geq 3$.*

*Proof.* Let $a_1, a_2, \ldots, a_k \in R$, where $k \geq 2$. Since $R$ is a ring satisfying the property ($\sharp$), there exists an idempotent $g \in R$ such that $a_k \in gRg$. We have

$$[a_1, a_2, \ldots, a_k]_k = [a_1, a_2, \ldots, a_k, g]_{k+1} \in [R, \ldots, R]_{k+1}.$$

Therefore, $[R, \ldots, R]_k \subseteq [R, \ldots, R]_{k+1}$ for all $k \geq 2$. In view of Lemma 3.2, we get $[R, \ldots, R]_n = \mathcal{I}([R, R])$ for all $n \geq 3$. □

Clearly, every unital ring has the property (♯). In view of [12, Lemma 2.4], every regular ring also satisfies the property (♯). The following is an immediate consequence of Theorem 3.7.

**Corollary 3.8.** *Let $R$ be either a unital ring or a regular ring. Then $[R, \ldots, R]_n = \mathcal{I}([R, R])$ for all $n \geq 3$.*

We next turn to answer Problem 1 for rings, which are generated by idempotents. Note that a ring $R$, which is generated by its idempotents, satisfies $R = R^2$. Given a nonempty subset $T$ of $R$, let $\overline{T}$ denote the subring of $R$ generated by $T$. We begin with a known result (see, for instance, [9, Fact 2]).

**Lemma 3.9.** *Let $R$ be a ring with an ideal $I$ and let $B$ be an additive subgroup of $R$. Then $[I, B] = [I, \overline{B}]$.*

**Theorem 3.10.** *Let $R$ be a ring generated by its idempotents. Then $[R, \ldots, R]_n = \mathcal{I}([R, R])$ for $n \geq 3$.*

*Proof.* Let $E$ denote the additive subgroup of $R$ generated by its all idempotents. By assumption, we have $\overline{E} = R$.

**Step 1**: $[R, R] \subseteq [R, \ldots, R]_{k+2}$ for $k \geq 1$. Let $x \in R$ and $e = e^2 \in R$. Then

$$[x, e] = [x, e_1, e_2, \ldots, e_{k+1}]_{k+2},$$

where $e_i = e$ for all $i$. This implies that $[x, e] \in [R, \ldots, R]_{k+2}$. That is, $[R, E] \subseteq [R, \ldots, R]_{k+2}$. It follows from Lemma 3.9 that $[R, E] = [R, \overline{E}] = [R, R]$. Therefore, $[R, R] \subseteq [R, \ldots, R]_{k+2}$.

**Step 2**: $[R, R]R \subseteq [R, \ldots, R]_{k+2}$ for $k \geq 1$. Let $w, x \in R$ and $e = e^2 \in R$. We claim that $[w, e]x \in [R, \ldots, R]_{k+2}$. We compute

$$\begin{aligned} [w, e]x &= [w, e, e]x \\ &= [w, e, ex] - [w, xe, e] + [wx, e, e] \\ &= [w, e_1, e_2, \ldots, e_k, ex]_{k+2} - [w, xe, e_1, e_2, \ldots, e_k]_{k+2} + [wx, e, e_1, e_2, \ldots, e_k]_{k+2} \\ &\in [R, \ldots, R]_{k+2}, \end{aligned}$$

where $e_i = e$ for all $i$. This implies that $[R, E]R \subseteq [R, \ldots, R]_{k+2}$. Recall that $[R, E] = [R, R]$. Therefore, $[R, R]R \subseteq [R, \ldots, R]_{k+2}$.

By Steps 1 and 2, $[R, R] + [R, R]R \subseteq [R, \ldots, R]_{k+2}$ for $k \geq 1$. By Lemma 2.2, $[R, R] + [R, R]R = \mathcal{I}([R, R])$ and, by (i) of Proposition 2.3, we have $[R, \ldots, R]_{k+2} = \mathcal{I}([R, R])$. □

**Example. (3)** Let $F$ be a field of characteristic 2, and let $V$ be the vector space over $F$ with an infinitely countable basis $v_1, v_2, \ldots, v_n, \ldots$. Let $R$ be the algebra generated by $V$ with multiplications: $v_i v_j = v_i$ for all $i, j \geq 1$. Then $R$ is generated by all its idempotents since $v_i^2 = v_i$ for all $i$. Also, $R = R^2$ and a direct computation shows that

$$\mathcal{I}([R,R]) = [R, \ldots, R]_k = \sum_{1 \leq i < j} F(v_i + v_j)$$

for all $k \geq 2$. Moreover, $\mathcal{I}([R,R])$ is nilpotent of index 2, $R\mathcal{I}([R,R]) = 0$ and $R/\mathcal{I}([R,R]) \cong F$. Indeed, the map $\phi\colon R \to F$ defined by $\phi(\sum_{i=1}^{s} \beta_i v_i) = \sum_{i=1}^{s} \beta_i$ where all $\beta_i \in F$ is an epimorphism with kernel $\mathcal{I}([R,R])$. In particular, $\mathcal{I}([R,R])$ is the unique maximal ideal of $R$.

Given $x, y \in R$, write $x = \sum_{i=1}^{s} \alpha_i v_i$ and $y = \sum_{j=1}^{t} \beta_j v_j$, where $\alpha_i, \beta_j \in F$ for all $i, j$. Then $xy = (\sum_{j=1}^{t} \beta_j)x$, implying that $[xy, x] = 0$. This means that $R$ is a PI-ring satisfying the identity $[XY, X]$. □

The following is well-known (see, for instance, [22, Lemma 2.1] with a short proof).

**Lemma 3.11.** *If $L$ is a Lie ideal of an arbitrary ring $R$, then $\mathcal{I}([L,L]) \subseteq L + L^2$.*

**Corollary 3.12.** *Let $R$ be a ring, and let $E$ be the additive subgroup of $R$ generated by all idempotents of $R$. Then*

$$[\overline{E}, \overline{E}] + R[\overline{E}, \overline{E}]R[\overline{E}, \overline{E}]R[\overline{E}, \overline{E}]R \subseteq [R, \ldots, R]_n$$

*for all $n \geq 3$.*

*Proof.* Let $n \geq 3$ be a positive integer. Since $\overline{E}$ is a ring generated by its idempotents, it follows from Theorem 3.10 that

$$\mathcal{I}_{\overline{E}}([\overline{E}, \overline{E}]) = [\overline{E}, \ldots, \overline{E}]_n \subseteq [R, \ldots, R]_n,$$

where $\mathcal{I}_{\overline{E}}([\overline{E}, \overline{E}])$ denotes the ideal of $\overline{E}$ generated by $[\overline{E}, \overline{E}]$. In view of Lemma 2.2,

$$\mathcal{I}_{\overline{E}}([\overline{E}, \overline{E}]) = [\overline{E}, \overline{E}] + \overline{E}[\overline{E}, \overline{E}] + \overline{E}[\overline{E}, \overline{E}]\overline{E}.$$

Note that $\overline{E}$ is both a subring and a Lie ideal of $R$. It follows from Lemma 3.11 that $R[\overline{E}, \overline{E}]R \subseteq \overline{E}$. Hence $[\overline{E}, \overline{E}] + R[\overline{E}, \overline{E}]R[\overline{E}, \overline{E}]R[\overline{E}, \overline{E}]R \subseteq [R, \ldots, R]_n$, as desired. □

A ring $R$ is called semiprime if, for $a \in R$, $aRa = 0$ implies $a = 0$. The semiprimeness of a ring $R$ is equivalent to saying that $R$ has no nonzero nilpotent one-sided ideal.

**Corollary 3.13.** *Let $R$ be a semiprime ring, and let $E$ be the additive subgroup of $R$ generated by all idempotents of $R$, and $n \geq 2$ a positive integer. Then either $ER \subseteq Z(R)$ or $[R, \ldots, R]_{2n}$ contains the nonzero ideal $R[\overline{E}, \overline{E}]R[\overline{E}, \overline{E}]R[\overline{E}, \overline{E}]R$.*

*Proof.* Assume that $[R, \ldots, R]_{2n}$ does not contain any nonzero ideal of $R$. Suppose first that $[E, E] \neq 0$. The semiprimeness of $R$ implies that

$$0 \neq (R[\overline{E}, \overline{E}])^4 \subseteq R[\overline{E}, \overline{E}]R[\overline{E}, \overline{E}]R[\overline{E}, \overline{E}]R \subseteq [R, \ldots, R]_{2n}.$$

Therefore, $R[\overline{E}, \overline{E}]R[\overline{E}, \overline{E}]R[\overline{E}, \overline{E}]R$ is a nonzero ideal of $R$ contained in $[R, \ldots, R]_{2n}$. This is a contradiction.

Therefore, $[E, E] = 0$. Let $e = e^2 \in R$ and $x \in R$. Then $e+ex(1-e)$ and $e+(1-e)xe$ are idempotents of $R$. Therefore, $[e, e + ex(1 - e)] = 0 = [e, e + (1 - e)xe]$, implying that $ex = xe$. That is, $e \in Z(R)$. Given $x, y, z \in R$ and $e = e^2 \in R$, we have

$$e[x, y, z] = [ex, y, z, e_1, \ldots, e_{2n-3}]_{2n} \in [R, \ldots, R]_{2n},$$

where $e_i = e$ for all $i$. This implies that $e[R, R, R] \subseteq [R, \ldots, R]_{2n}$. Clearly, $e[R, R, R]$ is an ideal of $R$ and hence $e[R, R, R] = 0$. In particular, $e[R, R, e] = 0$ and hence $e[R, R] = 0$. We get $eR \subseteq Z(R)$. This proves that $ER \subseteq Z(R)$, as desired. □

The following example constructed in the proof of [23, Theorem 1.2] shows that $[R, \ldots, R]_n = \mathcal{I}([R, R])$ holds for all $n \geq 2$ even when neither $R$ has the unity nor $R$ is regular.

**Example. (4)** Let $F$ be a field, and let $n \geq 3$ be a positive integer. Let

$$R := Fe_{21} + \sum_{2 \leq i \leq j \leq n} Fe_{ij} \subseteq \mathrm{M}_n(F).$$

Clearly, $R$ is a subring of $\mathrm{M}_n(F)$ and $R = R^2$. Therefore, $R$ is a PI-ring. Since $e_{21}R = 0$, $R$ is not a unital ring. A direct computation shows that

$$\mathcal{I}([R, R]) = [R, R] = Fe_{21} + \sum_{2 \leq i < j \leq n} Fe_{ij}.$$

Moreover, $[R, R] = [R, R, \ldots, R]_k$ for all $k \geq 3$, $\mathcal{I}([R, R])$ is nilpotent of index $n$, and $R/\mathcal{I}([R, R]) \cong F_1 \oplus \cdots \oplus F_{n-1}$, where $F_i = F$ for $i = 1, \ldots, n - 1$. □

# 4 $n$-Generalized commutator rings

As followed above (see [26]), by a commutator ring we mean a ring $R$ satisfying $R = [R, R]$. We study a generalization of commutator rings as follows.

**Definition 2.** Given a positive integer $n \geq 2$, a ring $R$ is called *an $n$-generalized commutator ring* provided $R = [R, \dots, R]_n$.

Therefore, a 2-generalized commutator ring just means a commutator ring (see [26]) and a 3-generalized commutator ring means a generalized commutator ring due to Herstein (see [14]). Clearly, the direct sum (resp. direct product) of $n$-generalized commutator rings is also an $n$-generalized commutator ring.

In view of Corollary 3.8, a unital ring $R$ is an $n$-generalized commutator ring where $n \geq 3$ if and only if $1 \in \mathcal{I}([R, R])$. In particular, every unital commutator ring is an $n$-generalized commutator ring for any $n \geq 3$. For commutator rings, see [26, Propositions 7, 12 and Theorem 13], and for algebras with a surjective inner derivation, see [30, Examples 1.1–1.6].

Given a unital ring $R$, let $A_n(R)$ denote the $n$-th Weyl algebra over $R$ (see [26, Definition 3]). Suppose that $R$ is a $\mathbb{Z}/p\mathbb{Z}$-algebra where $p$ is a prime integer, which is not a commutator ring. In view of [26, Proposition 8], $A_n(R)$ is not a commutator ring. However, since $1 \in [A_n(R), A_n(R)]$, it follows from Corollary 3.8 that $A_n(R)$ is an $m$-generalized commutator ring for any $m \geq 3$.

Let $R := \mathrm{M}_m(T)$, where $T$ is a unital ring and $m \geq 2$. In view of [21, Theorem 2.1], we have $R = \mathcal{I}([R, R])$. It follows from Corollary 3.8 that $R$ is an $n$-generalized commutator ring for all $n \geq 3$. Note that $[R, \dots, R]_n$ contains all commutators and generalized commutators if $n \geq 3$. Khurana and Lam proved a more strong result: Every element of $R$ is the sum of one commutator and a generalized commutator (see [17, (2) of Theorem A]). In addition, if $T$ is a PI-ring, then so is $R$ (see [28, Theorem 6.1.1]). In view of [2, Theorem 1], $R \neq [R, R]$ (that is, $R$ is a 3-generalized commutator ring but is not a commutator ring, see also [26, Proposition 6]).

**Examples. (5)** Let $A, B$ be unital rings, and let $M$ be a unital $(A, B)$-bimodule. Let $R := \mathrm{Tri}(A, M, B)$ be the triangular ring consisting of all elements $\begin{pmatrix} a & m \\ 0 & b \end{pmatrix}$ for $a \in A, m \in M, b \in B$ under the usual matrix operations. Note that $\begin{pmatrix} 0 & m \\ 0 & 0 \end{pmatrix} \in [R, R]$ since $[\begin{pmatrix} 1 & 0 \\ 0 & 0 \end{pmatrix}, \begin{pmatrix} 0 & m \\ 0 & 0 \end{pmatrix}] = \begin{pmatrix} 0 & m \\ 0 & 0 \end{pmatrix}$. Therefore, $R$ is a commutator ring if and only if both $A$ and $B$ are. Moreover,

$$\mathcal{I}([R,R]) = \left( \begin{array}{cc} \mathcal{I}([A,A]) & M \\ 0 & \mathcal{I}([B,B]) \end{array} \right).$$

Therefore, for a positive integer $n \geq 3$, $R$ is an $n$-generalized commutator ring if and only if both $A$ and $B$ are. □

**(6)** Given a positive integer $m$, let $T_m(R)$ denote the upper triangular subring of $\mathrm{M}_m(R)$, where $R$ is a unital ring. Then, for $m \geq 2$, we have $T_m(R) = \mathrm{Tri}(R, M, T_{m-1}(R))$, where $M = \{\sum_{j=2}^{m} a_{1j}e_{1j} \mid a_{1j} \in R\}$. In view of Corollary 3.8 and Example (4), applying the inductive argument we can prove that if $n \geq 2$, then $T_m(R)$ is an $n$-generalized commutator ring if and only if $R$ is. □

**(7)** Let $R := \mathbb{Z}\cdot I_m + \mathrm{M}_m(2\mathbb{Z})$ with $m > 1$, where $I_m$ is the identity matrix of $\mathrm{M}_m(\mathbb{Z})$. Then $R$ is a unital ring but $I_m \notin [R, \ldots, R]_n$ for $n > 1$, in particular, $[R, \ldots, R]_n \neq R$.

Indeed, it is clear that $I_m \notin [R, R]$. Suppose that $n \geq 3$. In view of Corollary 3.8, $[R, \ldots, R]_n = \mathcal{I}([R,R])$. Therefore, $\mathcal{I}([R,R]) \subseteq \mathrm{M}_m(4\mathbb{Z})$. In particular, $I_m \notin [R, \ldots, R]_n$, as desired. □

In a recent paper, Eroğlu prove that if $1 \in [R, R]$ then $R = \overline{[R, R]}$ (see [8, Theorem 1.1] and also [22, Theorem 1.3]). As an immediate consequence of Corollary 3.8, we have the following corollary.

**Corollary 4.1.** *Let $R$ be a unital ring, and $n \geq 3$ a positive integer. If $1 \in [R, \ldots, R]_n$, then $R$ is an $n$-generalized commutator ring.*

We continue to study Problem 1 in various ways.

**Proposition 4.2.** *Let $R$ be a ring, and $n \geq 2$ a positive integer. Then*

$$\big[K, R, \ldots, R\big]_n \subseteq [R, R, \ldots, R]_{n+1},$$

*where $K := [R, R^2] + R[R, R]$. Moreover, $K$ is an ideal of $R$. In addition, if $R = R^2$ then $K = \mathcal{I}([R, R])$.*

*Proof.* Given a positive integer $n \geq 2$, we let

$$A := \{a \in R \mid [a, R, \ldots, R]_n \subseteq [R, R, \ldots, R]_{n+1}\}.$$

Clearly, $A$ is an additive subgroup of $R$. We first claim that $x^2, x^3 \in A$ for all $x \in R$. Given $x \in R$, for $j = 1, 2$ we have

$$[x, x^j, y_1, \ldots, y_{n-1}]_{n+1} = [x^{j+1}, y_1, \ldots, y_{n-1}]_n$$

for all $y_1, \ldots, y_{n-1} \in R$. This proves that $x^2, x^3 \in A$.

Let $x, y, z \in R$. Then $xy + yx = (x+y)^2 - x^2 - y^2 \in A$. Therefore,

$$(xy)z + z(xy), (zx)y + y(zx) \in A$$

and so

$$[x, yz] = xyz - yzx = \big((xy)z + z(xy)\big) - \big((zx)y + y(zx)\big) \in A.$$

That is, $[R, R^2] \subseteq A$.

We next claim that $R[R, R] \subseteq A$. Since $(x^2y + yx^2) + (y^2x + xy^2) \in A$ and

$$(x+y)^3 - x^3 - y^3 = (x^2y + yx^2) + (y^2x + xy^2) + xyx + yxy \in A,$$

we get $xyx + yxy \in A$. Linearizing $xyx + yxy$ at $x$, we get $xyz + zyx \in A$. Then

$$x[y, z] = xyz - xzy = \big(xyz + zyx\big) - \big((zy)x + x(zy)\big) \in A.$$

That is, $R[R, R] \subseteq A$. Hence $[R, R^2] + R[R, R] \subseteq A$, as desired.

Recall that $K := [R, R^2] + R[R, R]$. Then

$$RK = R[R, R^2] + R^2[R, R] \subseteq R[R, R] \subseteq K.$$

In view of Lemma 2.2, $R[R, R]$ is an ideal of $R$. Therefore, $(R[R, R])R \subseteq R[R, R]$ and so

$$KR = [R, R^2]R + (R[R, R])R \subseteq [R^2, R^2] + R[R, R^2] + (R[R, R])R \subseteq K.$$

This proves that $K$ is an ideal of $R$. In addition, if $R = R^2$, then it follows from Lemma 2.2 that $K = [R, R] + R[R, R] = \mathcal{I}([R, R])$. □

**Corollary 4.3.** *Let $R$ be a noncommutative semiprime ring, $n \geq 2$ a positive integer. Then $\big[[R, R^2] + R[R, R], R, \ldots, R\big]_{2n-1}$ is a nonzero ideal of $R$, which is contained in $[R, R, \ldots, R]_{2n}$.*

*Proof.* Set $I := \big[[R, R^2] + R[R, R], R, \ldots, R\big]_{2n-1}$. Since $[R, R^2] + R[R, R]$ is an ideal of $R$ (see Proposition 4.2). In view of Eq.(1) and Eq.(2), $I$ is an ideal of $R$, which is contained in $[R, R, \ldots, R]_{2n}$. It suffices to claim that $I \neq 0$. Otherwise, we have $I = 0$.

Let $w, x, y, z \in R$. Then

$$0 = \big[w[x, y], z, z, \ldots, z\big]_{2n} = \big[w[x, y], z^{2n-1}\big]$$

for all $w, x, y, z \in R$. In view of [19, Theorem, p.19], $\big[w[x, y], z\big] = 0$ for all $w, x, y, z \in R$. That is, $\big[R[R, R], R\big] = 0$ and so $\big[R(R[R, R]), R\big] = 0$. This implies that $[R, R]R[R, R] = 0$. The semiprimeness of $R$ implies that $R$ is commutative, a contradiction. □

As a consequence of Proposition 4.2, the following characterizes $n$-generalized commutator rings for $n \geq 3$.

**Theorem 4.4.** *Let $R$ be a ring. Then the following are equivalent:*
*(i) $R = \mathcal{I}([R, R])$.*
*(ii) $R$ is an $n$-generalized commutator ring for all $n \geq 3$.*
*(iii) $R$ is an $n$-generalized commutator ring for some $n \geq 3$.*

*Proof.* "(i) $\Rightarrow$ (ii)": Assume that $R = \mathcal{I}([R, R])$. This implies that $R = R^2$ and hence $[R, R^2] = [R, R]$. By (ii) of Lemma 2.2, we have $[R, R^2] + R[R, R] = \mathcal{I}([R, R])$. In virtue of Proposition 4.2,

$$[R, R, \ldots, R]_{n-1} = \big[[R, R^2] + R[R, R], R, \ldots, R\big]_{n-1} \subseteq [R, R, \ldots, R]_n$$

for all $n \geq 3$. In view of Lemma 3.2, we get $R = \mathcal{I}([R, R]) = [R, R, \ldots, R]_n$, as desired.

It is trivial for "(ii) $\Rightarrow$ (iii)".

"(iii) $\Rightarrow$ (i)": Suppose that $R$ is an $n$-generalized commutator ring for some $n \geq 3$. In particular, by (i) of Proposition 2.3 we have $R = [R, R, \ldots, R]_n \subseteq \mathcal{I}([R, R])$ and so $R = \mathcal{I}([R, R])$. □

Since $R = \mathcal{I}([R, R])$ for any commutator ring $R$, the following corollary is an immediate consequence of Theorem 4.4.

**Corollary 4.5.** *Let $R$ be a commutator ring. Then $R$ is an $n$-generalized commutator ring for all $n \geq 3$.*

Taking into account what we have obtained so far, even a weaker version of the stated above Problem 1 is still unknown: Given a ring $R$ satisfying $R = R^2$, if $[R, R] = \mathcal{I}([R, R])$, is $[R, \ldots, R]_n$ equal to $\mathcal{I}([R, R])$ for any $n \geq 3$? A closely related question to the last corollary is of whether or not $R/I$ being a commutator ring for some nilpotent ideal $I$ of nilpotence index $k$ will imply that $R$ is an $n$-generalized commutator ring for some $n > 2$ depending on $k$? It is not in general true. Indeed, we choose a commutator ring $S$ and a nilpotent ring $I$ with nilpotence index 2 (that is, $I^2 = 0$ but $I \neq 0$). Let $R := S \oplus I$. Then $R/I \cong S$, implying that $R/I$ is a commutator ring. Suppose on the contrary that $R$ is an $n$-generalized commutator ring for some $n \geq 3$. Then $R = [R, \ldots, R]_n \subseteq R^n \subseteq R^2$, implying that $R = R^2$. However, $R^2 = S^2 \subsetneq S \oplus I = R$, a contradiction.

We answer the question with the necessary assumption that $R = R^2$ in the affirmative as a generalization of Corollary 4.5.

**Theorem 4.6.** *Let $R$ be a ring satisfying $R = R^2$. Suppose that $R/I$ is a commutator ring for some nilpotent ideal $I$ of $R$. Then $R$ is an $n$-generalized commutator ring for all $n \geq 3$.*

*Proof.* Since $I$ is a nilpotent ideal of $R$, $I^k = 0$ for some positive integer $k > 1$. By assumption, $R/I = [R/I, R/I]$. Therefore, $R = [R, R] + I$. Since $R = R^2$, we get $R = R^k$. This implies that

$$R = R^k = \big([R, R] + I\big)^k \subseteq \mathcal{I}([R, R]) + I^k = \mathcal{I}([R, R]).$$

That is, $R = \mathcal{I}([R, R])$. In view of Theorem 4.4, $R$ is an $n$-generalized commutator ring for all $n \geq 3$. □

Let $R$ be either a simple ring which is not a PI-ring or a commutator ring. Then $R/I$ is not a PI-ring for any proper ideal $I$ of $R$. The first case is clear. The latter case is then a consequence of [2, Theorem 1].

**Theorem 4.7.** *Let $R$ be a ring such that $R/I$ is not a PI-ring for any proper ideal $I$ of $R$. Then $R$ is an $n$-generalized commutator ring for all $n \geq 3$.*

*Proof.* In view of Lemma 2.2, $I := R[R, R]$ is an ideal of $R$. Suppose that $I$ is a proper ideal of $R$. Then $R/I$ satisfies the polynomial identity $X_1[X_2, X_3]$, a contradiction. Therefore, $R[R, R] = R$, implying that $R = \mathcal{I}([R, R])$. In view of Theorem 4.4, $R$ is an $n$-generalized commutator ring for all $n \geq 3$. □

It is worthwhile noticing that the so-constructed rings in Examples (3) and (4) are PI-rings that, by virtue of Theorem 4.4, are not $n$-generalized commutator rings for any $n > 1$, because $R$ is not equal to $\mathcal{I}([R, R])$.

The following lemma will play a key role in the sequel.

**Lemma 4.8.** *Let $R$ be a ring and $n \geq 3$ a positive integer. Then the following hold:*
*(i)* $[z^{n-1}, R] \subseteq [R, \ldots, R]_n$ *for all* $z \in R$.
*(ii)* $z^{n-1}[R, R] \subseteq [R, \ldots, R]_n$ *for all* $z \in R$.
*(iii)* $\sum_{z \in R} R z^{n-1}[z^{n-1}, R]R \subseteq [R, \ldots, R]_n$.

*Proof.* Let $L := [R, \ldots, R]_n$, where $n \geq 3$. It is known that $L$ is a Lie ideal of $R$. Given $x, y, z \in R$, we have

$$[x, z^{n-1}] = [x, z_1, \ldots, z_{n-2}, z_{n-1}]_n \in L, \tag{3}$$

where $z_i = z$ for all $i$. Moreover,

$$[x, y, z^{n-1}] = xyz^{n-1} - z^{n-1}yx = [x, y, z_1, \ldots, z_{n-4}, z_{n-3}, z^2]_n \in L, \tag{4}$$

where $z_i = z$ for all $i$. It follows from Eq.(3) and Eq.(4) that

$$z^{n-1}[x, y] = [z^{n-1}, xy] + [x, y, z^{n-1}] \in L.$$

Up to now, we have proved that $[R, z^{n-1}] \subseteq L$ and $z^{n-1}[R, R] \subseteq L$ for all $z \in R$. Therefore, we have established (i) and (ii).

We now prove (iii). Given $z \in R$, applying (i) and (ii) we have

$$z^{n-1}[z^{n-1}, R]R \subseteq z^{n-1}[z^{n-1}R, R] + (z^2)^{n-1}[R, R] \subseteq L$$

and hence

$$Rz^{n-1}[z^{n-1}, R]R \subseteq \big[R, z^{n-1}[z^{n-1}, R]R\big] + z^{n-1}[z^{n-1}, R]R^2 \subseteq [R, L] + L \subseteq L.$$

Therefore, $\sum_{z \in R} Rz^{n-1}[z^{n-1}, R]R \subseteq [R, \ldots, R]_n$, as desired. □

Let $R$ be a noncommutative simple ring, and $k$ a fixed positive integer. Applying [5, Theorem 2] (i.e. Theorem 7.5 below), we can show that the subring of $R$ generated by all elements $z^k$ for $z \in R$ is equal to the whole ring $R$. The following theorem answers Problem 1 affirmatively for rings of such type.

**Theorem 4.9.** *Let $R$ be a ring, $n \geq 3$ a positive integer. Suppose that the ring $R$ is generated by all elements $z^{n-1}$ for $z \in R$. Then $[R, \ldots, R]_n = \mathcal{I}([R, R])$.*

*Proof.* Set $L := [R, \ldots, R]_n$. We let $A$ denote the additive subgroup of $R$ generated by all elements $z^{n-1}$ for $z \in R$. By assumption, we have $\overline{A} = R$. By (i) of Lemma 4.8, we get $[A, R] \subseteq L$. It follows from Lemma 3.9 that $[R, R] = [\overline{A}, R] = [A, R] \subseteq L$.

Let $z \in R$. In view of (ii) of Lemma 4.8, we have $z^{n-1}[R, R] \subseteq L$. Therefore,

$$[z^{n-1}, R]R \subseteq z^{n-1}[R, R] + [z^{n-1}R, R] \subseteq L + [R, R] = L.$$

That is, $[A, R]R \subseteq L$ and so $[R, R]R = [\overline{A}, R]R = [A, R]R \subseteq L$. By (ii) of Lemma 2.2 and (i) of Proposition 2.3, we have

$$\mathcal{I}([R, R]) = [R, R] + [R, R]R \subseteq L \subseteq \mathcal{I}([R, R]).$$

This proves that $\mathcal{I}([R, R]) = L$, as desired. □

**Example. (8)** Let $R := \mathrm{M}_n(2\mathbb{Z})$. Any nonzero ideal of $R$ where $n \geq 2$ is not a $k$-generalized commutator ring for $k \geq 2$. Moreover, $R \supsetneq 2R \supsetneq 2^2R \supsetneq \cdots$ is an infinite descending chain of ideals of $R$.

Indeed, let $N$ be a nonzero ideal of $\mathrm{M}_n(2\mathbb{Z})$. Choose a nonzero element $x := \sum_{1\le i,j\le n} a_{ij}e_{ij} \in N$, where $a_{ij} \in \mathbb{Z}$ for all $i, j$. Clearly, there exists a positive integer $m$ such that $x \in \mathrm{M}_n(2^m\mathbb{Z})$ but $x \notin \mathrm{M}_n(2^{m+1}\mathbb{Z})$. Suppose on the contrary that $N$ is a $k$-generalized commutator ring for some $k \ge 2$. In particular, $N = N^2$ and so $N = N^{m+1}$. This implies that $N \subseteq \mathrm{M}_n(2^{m+1}\mathbb{Z})$ and so $x \in \mathrm{M}_n(2^{m+1}\mathbb{Z})$, a contradiction. The final assertion is then clear. □

Motivated by Example (8), a ring $R$ is said to satisfy the descending chain condition (d.c.c. for short) on ideals if every nonempty set of ideals of $R$ contains a minimal element. It is equivalent to saying that each descending chain of ideals $I_1 \supseteq I_2 \supseteq I_3 \supseteq \cdots$ must be stationary.

Before proving our next theorem, we need the following technical lemma.

**Lemma 4.10.** *Let $R$ be a semiprime ring with a right ideal $\rho$. If $[a, \rho] \subseteq Z(R)$ where $a \in R$, then $\rho[a, R] = 0$.*

*Proof.* Since $[a, \rho]a = [a, \rho a] \subseteq [a, \rho] \subseteq Z(R)$, we get $\big[[a, \rho]a, R\big] = 0$. This implies that $[a, \rho][a, R] = 0$. In particular, $[a, \rho][a, R^2] = 0$, implying that $[a, \rho]R[a, R] = 0$. Therefore, $[a, \rho]R[a, \rho] = 0$. The semiprimeness of $R$ implies that $[a, \rho] = 0$ and so $[a, \rho R] = 0$. We thus get $\rho[a, R] = 0$, as desired. □

**Theorem 4.11.** *Let $R$ be a noncommutative semiprime ring, and let*

$$K_I := \sum_{z\in I} Iz^{n-1}[z^{n-1}, I]I$$

*for an ideal $I$ of $R$, where $n \ge 3$ is a positive integer. Then the following hold:*

*(i) Given an ideal $I$ of $R$, $K_I$ is an ideal of $R$ and $K_I \subseteq [I, \ldots, I]_n$, and if $[I, I] \ne 0$ then $[K_I, K_I] \ne 0$.*

*(ii) If $R$ satisfies the d.c.c. on ideals, then there exists a nonzero ideal $N$ of $R$ such that $N$ is an $n$-generalized commutator ring for all $n \ge 3$.*

*Proof.* (i) Fix a positive integer $n \ge 3$. Denote by $Q_I$ the Martindale symmetric ring of quotients of $I$ (see [1] for its definition). Since $R$ is a semiprime ring, so is $I$. It follows from the proof of Lemma 4.8 that

$$K_I = \sum_{z\in I} Iz^{n-1}[z^{n-1}, I]I \subseteq [I, \ldots, I]_n.$$

Clearly, $K_I$ is also an ideal of $R$. We claim that if $I$ is noncommutative, so is $K_I$. Otherwise, we have $[K_I, K_I] = 0$. Since $K_I$ is an ideal of $I$, this implies that $K_I \subseteq Z(I)$. Therefore,

$$\big[xz^{n-1}[z^{n-1}, y]v, w\big] = 0 \tag{5}$$

for all $v, w, x, y, z \in I$. In view of [1, Theorem 6.4.1], $I$ and $Q_I$ satisfy the same GPIs. Hence Eq.(5) holds for all $v, w, x, y, z \in Q_I$. Replacing $x, v$ by 1 in Eq.(5), we get $z^{n-1}[z^{n-1}, y] \in Z(I)$ for all $y, z \in I$. That is, $\big[z^{n-1}, z^{n-1}I\big] \subseteq Z(I)$ for all $z \in I$. In view of Lemma 4.10, $z^{n-1}I[z^{n-1}, R] = 0$ for all $z \in I$. Hence $[z^{n-1}, I]I[z^{n-1}, I] = 0$ for all $z \in I$. The semiprimeness of $I$ asserts that $[z^{n-1}, I] = 0$ for all $z \in I$. In view of [19, Theorem, p.19], $[z, I] = 0$ for all $z \in I$. Therefore, $I$ is commutative, a contradiction. This proves (i).

(ii) We let

$$I_0 := K_R \ \text{ and } \ I_n := K_{I_{n-1}} \ \text{ for } \ n = 1, 2, \ldots.$$

Then $I_j$ is an ideal of $R$ for all $j \geq 0$ and $I_0 \supseteq I_1 \supseteq I_2 \supseteq I_3 \supseteq \cdots$. Moreover, by $[R, R] \neq 0$, we have $[I_j, I_j] \neq 0$ for all $j \geq 0$. Since $R$ satisfies the d.c.c. on ideals, there exists a positive integer $k$ such that $I_k = I_s$ for all $s > k$. We let $N := I_k$. Then

$$0 \neq N = I_k = I_{k+1} = K_{I_k} = \sum_{z \in I_k} I_k z^{n-1}[z^{n-1}, I_k]I_k \subseteq [I_k, \ldots, I_k]_n = [N, \ldots, N]_n,$$

implying that $N = [N, \ldots, N]_n$, as desired. In view of Theorem 4.4, $N = [N, \ldots, N]_n$ for all $n \geq 3$. □

**Corollary 4.12.** *Let $R$ be a noncommutative semiprime ring, $n \geq 3$ a positive integer. Then $[R, \ldots, R]_n$ contains a nonzero ideal $W$ of $R$ such that $[W, W] \neq 0$.*

*Proof.* We let $W := \sum_{z \in R} Rz^{n-1}[z^{n-1}, R]R$. In view of (i) of Theorem 4.11, $W$ is a nonzero ideal of $R$ and $[W, W] \neq 0$, as desired. □

Clearly, every noncommutative simple ring satisfies the d.c.c. on ideals. Applying (ii) of Theorem 4.11, we have the following corollary, which is a generalization of [14, Theorem 4] for $n = 3$. Of course, it is also a consequence of Corollary 4.12.

**Corollary 4.13.** *Let $R$ be a noncommutative simple ring. Then $R$ is an $n$-generalized commutator ring for all $n \geq 3$.*

Applying Theorem 4.4, we give an alternative proof for (ii) of Theorem 4.11. Indeed, let $\Sigma := \{I \lhd R \mid [I, I] \neq 0\}$, where by $I \lhd R$ we mean that $I$ is an ideal of $R$. By assumption $[R, R] \neq 0$, we have $R \in \Sigma$. Since $R$ satisfies the d.c.c. on ideals, there exists a minimal element, say $N$, in $\Sigma$. Then $[N, N] \neq 0$. The semiprimeness of $R$ implies that $\big[[N, N], [N, N]\big] \neq 0$ (see Lemma 7.3). By the semiprimeness of the ring $N$ again, we have $\big[N[N, N]N, N[N, N]N\big] \neq 0$. Since $N[N, N]N$ is also an ideal of $R$, it follows that $N[N, N]N \in \Sigma$. The minimality of $N$ in $\Sigma$ implies that $N = N[N, N]N$ and so $N = \mathcal{I}_N([N, N])$, the ideal of $N$ generated by $[N, N]$. In view of Theorem 4.4, $N$ is an $n$-generalized commutator ring for all $n \geq 3$, as desired.

We end this section with an example. In view of Corollary 3.8 and Theorem 3.7, the existence of idempotents in the considered rings seems to play an important role for answering Problem 1 affirmatively. The following, however, shows that the existence of idempotents is not essential to Problem 1.

A ring $R$ is said to satisfy the ascending chain condition (a.c.c. for short) on ideals if every nonempty set of ideals of $R$ contains a maximal element. It is equivalent to saying that each ascending chain of ideals $I_1 \subseteq I_2 \subseteq I_3 \subseteq \cdots$ must be stationary.

**Example. (9)** There exists a nil PI-ring $R$ such that $R = R^2$, $R \neq [R, R]$ and $R$ is a $k$-generalized commutator ring for all $k \geq 3$. Moreover, neither $R$ satisfies the a.c.c. on ideals nor $R$ satisfies the d.c.c. on ideals.

Let $T$ be the commutative algebra over a field F with the symbols $v_\alpha$'s, where $0 < \alpha < 1$, as a basis. The multiplications of these elements $v_\alpha$'s for $0 < \alpha < 1$ are defined by $v_\alpha v_\beta = v_{\alpha+\beta}$ if $\alpha + \beta < 1$ and $v_\alpha v_\beta = 0$ if $\alpha + \beta \geq 1$. Clearly, $T$ is a nil commutative algebra over $F$. Let $R := \mathrm{M}_n(T)$, where $n \geq 2$. Note that, for any $x \in R$, there exists $0 < \alpha < 1$ such that $x \in v_\alpha R$. Therefore, $x \in R^2$ and $x^m = 0$ if $m\alpha \geq 1$. This proves that $R = R^2$ and $R$ is a nil ring. Since $T$ is commutative, $R$ is a PI-ring (see [28, Theorem 6.1.1]). In view of [2, Theorem 1], $R \neq [R, R]$.

Let $a_1, \ldots, a_k \in R$, where $k \geq 2$. There exists $0 < \alpha < 1$ such that $a_k \in v_\alpha R$. Write $a_k = v_\alpha a'_k$, where $a'_k \in R$. We have

$$[a_1, \ldots, a_k]_k = [a_1, \ldots, a_{k-1}, v_{\frac{\alpha}{2}} a'_k, v_{\frac{\alpha}{2}}]_{k+1} \in [R, \ldots, R]_{k+1}.$$

Therefore, $[R, \ldots, R]_k \subseteq [R, \ldots, R]_{k+1}$ for all $k \geq 2$. In view of Lemma 3.2, we get $[R, \ldots, R]_k = \mathcal{I}([R, R])$ for all $k \geq 3$. To prove that $R$ is a $k$-generalized commutator ring for all $k \geq 3$, it suffices to claim that $R = \mathcal{I}([R, R])$.

Let $0 < \alpha < 1$ and $1 \leq i, j \leq n$ with $i \neq j$. Then

$$v_\alpha e_{ij} = [v_{\frac{\alpha}{2}} e_{ii}, v_{\frac{\alpha}{2}} e_{ij}] \in [R, R].$$

On the other hand, we have $v_\alpha e_{ii} = [v_{\frac{\alpha}{3}} e_{ii}, v_{\frac{\alpha}{3}} e_{ij}] v_{\frac{\alpha}{3}} e_{ji} \in [R, R]R$. Hence $R = [R, R] + [R, R]R = \mathcal{I}([R, R])$, as desired.

Finally, we have

$$\mathcal{I}(\{v_{\frac{1}{2}}\}) \subsetneq \mathcal{I}(\{v_{\frac{1}{2^2}}\}) \subsetneq \mathcal{I}(\{v_{\frac{1}{2^3}}\}) \subsetneq \cdots$$

and

$$\mathcal{I}(\{v_{1-\frac{1}{2}}\}) \supsetneq \mathcal{I}(\{v_{1-\frac{1}{2^2}}\}) \supsetneq \mathcal{I}(\{v_{1-\frac{1}{2^3}}\}) \supsetneq \cdots$$

This proves that neither $R$ satisfies the a.c.c. on ideals nor $R$ satisfies the d.c.c. on ideals. □

# 5 $n$-Generalized Lie ideals

Let $n \geq 3$ be a positive integer. We have proved that if $R$ is a noncommutative prime ring, $[R, \ldots, R]_n$ contains a nonzero ideal of $R$ (see Corollary 4.12). In particular, if $R$ is a noncommutative simple ring, then $R$ is an $n$-generalized commutator ring (see Corollary 4.13).

Let $L := [R, \ldots, R]_n$. Then $[L, R, \ldots, R]_n \subseteq [R, \ldots, R]_n = L$. We will study these results above from this viewpoint. Our present study is also motivated by [15] and [18]. Precisely, in 1955 Herstein determined the Lie structure of simple rings (see [15]). In 1972 Lanski and Montgomery extended Herstein's theorem to the context of prime rings (see [18, Theorem 13] and the references therein).

We define $n$-generalized Lie ideals of rings for $n \geq 2$, which will coincide with Lie ideals if $n = 2$.

**Definition 3.** By an $n$-generalized Lie ideal of a ring $R$ (at the $(r+1)$-th position with $r \geq 0$) we mean an additive subgroup $A$ of $R$ satisfying $[x_1, \ldots, x_r, a, y_1, \ldots, y_s]_n \in A$ for all $x_i, y_j \in R$ and all $a \in A$, where $r + s = n - 1$.

Clearly, every ideal of $R$ is an $n$-generalized Lie ideal of $R$. Moreover, $[R, \ldots, R]_n$ is also an $n$-generalized Lie ideal of $R$. Note $[R, R]$ is a Lie ideal of $R$ but it does not in general contain a nonzero ideal of $R$. We are now ready to state the main theorem.

**Theorem 5.1.** *Let $R$ be a noncommutative prime ring and $n \geq 3$ a positive integer. Then every nonzero $n$-generalized Lie ideal of $R$ contains a nonzero ideal.*

The following assertion gives a generalization of Corollary 4.12 for prime rings.

**Corollary 5.2.** *Let $R$ be a noncommutative prime ring and $n \geq 3$ a positive integer. If $A$ is a nonzero $n$-generalized Lie ideal of $R$ (at the $(r+1)$-th position with $r \geq 0$), so is $[R_1, \ldots, R_r, A, R_1, \ldots, R_s]_n$ where $R_i = R$ for all $i$.*

*Proof.* We let $K := [R_1, \ldots, R_r, A, R_1, \ldots, R_s]_n$. Since $A$ is a nonzero $n$-generalized Lie ideal of $R$ (at the $(r+1)$-th position with $r \geq 0$), we have $K \subseteq A$. Therefore,

$$[R_1, \ldots, R_r, K, R_1, \ldots, R_s]_n \subseteq [R_1, \ldots, R_r, A, R_1, \ldots, R_s]_n = K,$$

implying that $K$ is a nonzero $n$-generalized Lie ideal of $R$ (at the $(r+1)$-th position with $r \geq 0$). In view of Theorem 5.1, $A$ contains a nonzero ideal, say $I$, of $R$. Then $[I, \ldots, I]_n \subseteq K$. The primeness of $R$ implies that $I$ is also a prime ring. Since $R$ is noncommutative, so is the prime ring $I$. In view of Theorem 5.1, we get $[I, \ldots, I]_n \neq 0$ and hence $K \neq 0$, as desired. □

The proof of Theorem 5.1 will be given in the next two sections. Some related questions on generalized commutators and their relationship with noncommutative polynomials are also discussed in the final section.

# 6 A special case: $[R, A, R] \subseteq A$

Throughout this section and the next one, unless specially stated, $R$ always denotes a prime ring with extended $C$, and let $Q$ be the Martindale symmetric ring of quotients of $R$. Recall that $Q$ is also a prime ring and that $C$ is a field. We refer the reader to the book [1] for details.

The aim of this section is to prove the following proposition, which is a special case of Theorem 5.1.

**Proposition 6.1.** *Let $R$ be a noncommutative prime ring and let $A$ be a nonzero additive subgroup of $R$. If $[R, A, R] \subseteq A$ then $A$ contains a nonzero ideal of $R$.*

We begin with the following lemma (see [24, Lemma 2.10] and [7, Lemma 3]).

**Lemma 6.2.** *Let $A$ be an additive subgroup of an arbitrary ring $R$. If $[R, A, R] \subseteq A$, then*

$$\sum_{a,\, b\in A,\, x\in R} R(axb - bxa)R \subseteq A.$$

*Proof.* Let $a, b \in A$ and $x, y, z \in R$. Then

$$\begin{aligned} & y(axb - bxa)z \\ &= (yax)bz - zb(yax) + (zby)ax - xa(zby) + (xaz)by - yb(xaz) \\ &= [yax, b, z] + [zby, a, x] + [xaz, b, y] \in A. \end{aligned}$$

□

Let $\mathbf{X} = \{X_1, X_2, \cdots\}$ be an infinitely countable set. We denote by $C\{\mathbf{X}\}$ the free $C$-algebra in noncommutative indeterminates in the set $\mathbf{X}$. We let $Q *_C C\{\mathbf{X}\}$ stand for the free product over $C$ of $C$-algebras $Q$ and $C\{\mathbf{X}\}$. By a generalized polynomial (GP for short) with coefficients in $Q$ we just mean an element in $Q *_C C\{\mathbf{X}\}$. For $f(X_1, \ldots, X_n) \in Q *_C C\{\mathbf{X}\}$, we just mean that $X_1, \ldots, X_n$ are the only indeterminates occurring in $f$. A polynomial $f(X_1, \ldots, X_n) \in Q *_C C\{\mathbf{X}\}$ is called a generalized polynomial identity (GPI for short) for $A$, an additive subgroup of $Q$, if $f(x_1, \ldots, x_n) = 0$ for all $x_i \in A$.

In order to simplify the proof below, we need the following lemma.

**Lemma 6.3.** *([20, Lemma 2.3]) Let* $h(X_1, \ldots, X_t) \in Q *_C C\{\mathbf{X}\} \setminus \{0\}$*. For a nonzero ideal* $I$ *of* $R$ *we let*

$$\mathcal{A} = \{h(x_1, \ldots, x_t) \mid x_i \in Q\} \text{ and } \mathcal{A}_I = \{h(x_1, \ldots, x_t) \mid x_i \in I\}.$$

*Then* $\dim_C \mathcal{A}C < \infty$ *if and only if* $\dim_C \mathcal{A}_I C < \infty$*. In this case,* $\mathcal{A}C = \mathcal{A}_I C$*.*

**Lemma 6.4.** *If* $xay - yax \in Ca$ *for all* $x, y \in R$*, where* $0 \neq a \in R$*, then* $R$ *is commutative.*

*Proof.* In view of Lemma 6.3, $xay - yax \in Ca$ for all $x, y \in Q$. By taking $x = 1$, we have $ay - ya \in Ca$ for all $y \in Q$. Therefore, $[a, [a, y]] = 0$ for all $y \in Q$.

**Case 1**: $\operatorname{char} R \neq 2$. It follows from [27, Theorem 1] that $a \in C$. Then $xay - yax = a[x, y] \in Ca = C$ for all $x, y \in Q$. Therefore, $[x, y] \in C$ for all $x, y \in Q$. This implies that $Q$ is commutative. In particular, $R$ is commutative.

**Case 2**: $\operatorname{char} R = 2$. Then $a^2 \in C$. If $a \in C$, then we are done as given in Case 1. Assume that $a \notin C$. Choose an element $y \in Q$ such that $[a, y] \neq 0$. Therefore, $[a, y] = \beta a$ for some nonzero $\beta \in C$. Then $a[a, y] = [a, ay] \in Ca$. That is, $\beta a^2 = \alpha a$ for some $\alpha \in C$. Since $a^2 \in C$ but $a \notin C$, we get $\alpha = 0$ and so $a^2 = 0$. So $a(ay - ya) \in Ca^2 = \{0\}$ for all $y \in Q$, implying that $aya = 0$ for all $y \in Q$. The primeness of $Q$ implies that $a = 0$, a contradiction. □

We are now ready to prove the main proposition in this section.

*Proof of Proposition 6.1.* Suppose that $[R, A, R] \subseteq A$. In view of Lemma 6.2,

$$\sum_{a,b \in A, x \in R} R(axb - bxa)R \subseteq A.$$

Therefore, either $A$ contains a nonzero ideal of $R$ or $axb = bxa$ for all $a, b \in A$ and all $x \in R$. It suffices to consider the latter case. By [25, Theorem 1], there exists $0 \neq a \in A$ such that $A \subseteq Ca$. In particular, $[x, a, y] = xay - yax \in Ca$ for all $x, y \in R$. It follows from Lemma 6.4 that $R$ is commutative, a contradiction. □

# 7 Proofs of Theorem 5.1

Recall that, unless specially stated, $R$ always denotes a prime ring with extended centroid $C$.

**Lemma 7.1.** *If* $R$ *is noncommutative and* $0 \neq a \in R$*, then* $[R, x]xa \neq 0$ *for some* $x \in R$*.*

*Proof.* Otherwise, $[R, x]xa = 0$ for all $x \in R$. Let $y, z, x \in R$. Then $[yz, x]xa = 0$. Since $y[z, x]xa = 0$, we have $[y, x]zxa = 0$. That is, $[R, x]Rxa = 0$. The primeness of $R$ implies that $[R, x] = 0$ or $xa = 0$. Hence $R$ is the union of its two additive subgroups $Z(R)$ and $\{x \in R \mid xa = 0\}$. Since $R$ is not commutative, $R = \{x \in R \mid xa = 0\}$ and so $Ra = 0$. Therefore, $a = 0$ follows, a contradiction. □

**Lemma 7.2.** *If $[w, x]x \in Z(R)$ for all $w, x \in R$, then $R$ is commutative.*

*Proof.* Since $R$ and $Q$ satisfy the same GPIs (see [1, Theorem 6.4.1]), $[w, x]x \in C$ for all $w, x \in Q$. Replacing $x$ by $x + 1$, we get $[w, x + 1](x + 1) \in C$ for all $w, x \in Q$. Therefore, $[w, x] \in C$ for all $w, x \in Q$, implying that $Q$ is commutative. In particular, $R$ is commutative. □

The following is well-known, but is listed here only for completeness of the exposition.

**Lemma 7.3.** *Let $I$ be a nonzero ideal of a semiprime ring $R$. Suppose that $[a, [I, I]] = 0$, where $a \in R$. Then $[a, I] = 0$. In addition, if either $a \in I$ or $R$ is a prime ring, then $a \in Z(R)$.*

*Proof.* Let $x, y \in I$. Then $xa \in I$ and so

$$0 = \big[a, [xa, y]\big] = \big[a, x[a, y] + [x, y]a\big] = \big[a, x[a, y]\big].$$

Therefore, $\big[a, xR[a, x]\big] = 0$ and so $[a, x]R[a, x] = 0$. The semiprimeness of $R$ implies that $[a, x] = 0$. That is, $[a, I] = 0$, as desired. It is clear that $a \in Z(R)$ if $a \in I$. □

**Lemma 7.4.** *Let $L_1$ and $L_2$ be two proper Lie ideals of $R$. If $R$ is noncommutative, then $[L_1, L_2] \neq 0$.*

*Proof.* Since the intersection of two nonzero ideals of $R$ is also a nonzero ideal, we can choose a nonzero ideal $I$ of $R$ such that $[I, I] \subseteq L_i$ for $i = 1, 2$. Suppose on the contrary that $[L_1, L_2] = 0$. Then $\big[[I, I], [I, I]\big] = 0$, implying $[I, I] \subseteq Z(R)$ by Lemma 7.3. This implies that $[R, [I, I]] = 0$ and hence $R \subseteq Z(R)$. That is, $R$ is commutative, a contradiction. □

Let $f(X_1, \ldots, X_m)$ be a multilinear polynomial over $\mathbb{Z}$, the ring of integers, in noncommuting indeterminates $X_1, \ldots, X_m$. Then

$$\big[y, f(x_1, \ldots, x_m)\big] = \sum_{k=1}^{m} f(x_1, \ldots, [y, x_k], \ldots, x_m)$$

for all $x_i, y \in R$. Therefore, the additive subgroup of $R$ generated by all elements $f(x_1, \ldots, x_m)$ for $x_i \in R$ is a Lie ideal of $R$. We will use the basic fact in the proof below.

For $f \in C\{X_1, \ldots, X_n\}$ and $I$ an ideal of $R$, let $\mathrm{Add}_f(I)$ denote the additive subgroup of $RC$ generated by all elements $f(x_1, \ldots, x_n)$ for $x_i \in I$. We say that $f$ is central-valued on $RC$ if $f(x_1, \ldots, x_n) \in C$ for all $x_i \in RC$. We also need the following theorem in our proof.

**Theorem 7.5.** *([5, Theorem 2]) Let $R$ be a prime ring and $I$ a nonzero ideal of $R$. Suppose that $f(X_1, \ldots, X_n) \in C\{X_1, \ldots, X_n\}$, which is not central-valued on $RC$. Then $[M, R] \subseteq \mathrm{Add}_f(I)$ for some nonzero ideal $M$ of $R$ except when $R \cong \mathrm{M}_2(\mathrm{GF}(2))$ and $\mathrm{Add}_f(R) = \{0, e_{12} + e_{21}, 1 + e_{12}, 1 + e_{21}\}$ or $\{0, 1, e_{11} + e_{12} + e_{21}, e_{22} + e_{12} + e_{21}\}$.*

We are now ready to prove Theorem 5.1.

*Proof of Theorem 5.1.* Let $A$ be an $n$-generalized Lie ideal of $R$ at the $(r+1)$-th position, where $0 \leq r \leq n$. That is,

$$\begin{aligned} &[x_1, \ldots, x_r, a, y_1, \ldots, y_s]_n \\ &= (x_1 x_2 \cdots x_r) a (y_1 \cdots y_s) - (y_s \cdots y_1) a (x_r x_{r-1} \cdots x_1) \in A \end{aligned} \tag{6}$$

for all $x_i, y_j \in R$ and all $a \in A$, where $r + s = n - 1$. In view of Proposition 6.1, we are done if $n = 3$ and $r = s = 1$. Therefore, we can always assume that either $r > 1$ or $s > 1$. By symmetry, we may assume that $r > 1$. In view of Eq.(6), we have

$$\begin{aligned} &[x_1 w, x_2 \ldots, x_r, a, y_1, \ldots, y_s]_n - [x_1, w x_2, \ldots, x_r, a, y_1, \ldots, y_s]_n \\ &= y_s \cdots y_1 a x_r x_{r-1} \cdots x_3 [w, x_2 x_1] \in A \end{aligned} \tag{7}$$

for all $x_i, y_j, w \in R$ and all $a \in A$.

**Case 1**: $s \geq 1$ and $r \geq 3$. In this case, let $L$ be the additive subgroup of $R$ generated by all elements $x_r x_{r-1} \cdots x_3 [w, x_2 x_1]$ for $x_1, \ldots, x_r, w \in R$. Recall that $R$ is a noncommutative prime ring. Clearly, $L$ is both a Lie ideal and a nonzero left ideal of $R$. In view of Lemma 2.2, $RL$ is a nonzero ideal of $R$ and $RL \subseteq L$. By Eq.(7), we have $0 \neq R^s ARL \subseteq A$, where $R^s ARL$ is an ideal of $R$, as desired.

**Case 2**: $s \geq 1$ and $r = 2$. In view of Eq.(7), we get $y_s \cdots y_1 a [w, x_2 x_1] \in A$ for all $y_j, x_1, x_2, w \in R$. That is, $\lambda := R^s A [R, R^2] \subseteq A$. By Eq.(6), we have

$$x_1 x_2 a (y_1 \cdots y_s) - (y_s \cdots y_1) a x_2 x_1 \in A \tag{8}$$

for all $x_1, x_2, y_j \in R$ and all $a \in A$. Choose $x_1 \in \lambda \subseteq A$. Then $(y_s \cdots y_1) a x_2 x_1 \in A$ and so, by Eq.(8), $x_1 x_2 a (y_1 \cdots y_s) \in A$. This implies that $\lambda RAR^s \subseteq A$. Clearly, $\lambda RAR^s$is a nonzero ideal of $R$, as desired.

**Case 3**: $s = 0$. Suppose first that $r \geq 3$. We have

$$(x_1 x_2 \cdots x_r)a - a(x_r x_{r-1} \cdots x_1) \in A$$

for all $x_1, x_2, \ldots, x_r \in R$ and all $a \in A$. We keep $L$ as given in Case 1. Then $L$ contains a nonzero ideal, say $M$, of $R$. Moreover, $0 \neq AM \subseteq A$. Choose $x_1 \in M$ and by the fact that $s = 0$, we get $a(x_r x_{r-1} \cdots x_1) \in AM \subseteq A$ and hence $MR^{r-1}A \subseteq A$. Then

$$MR^{r-1}AM \subseteq AM \subseteq A.$$

Clearly, $MR^{r-1}AM$ is a nonzero ideal of $R$.

Suppose next that $r = 2$. Then

$$[x_1 w, x_2, a] - [x_1, w x_2, a] = a[w, x_2 x_1] \in A$$

for all $x_1, x_2 \in R$ and all $a \in A$. Therefore, $AL_1 \subseteq A$, where $L_1 := [R, R^2]$. Let $w, x_2 \in R$ and $a \in A$. We also compute

$$\big[[w, x_2], x_2, a\big] = [w, x_2]x_2 a - a x_2 [w, x_2] \in A.$$

Note that $ax_2[w, x_2] = a[x_2 w, x_2] \in A[R^2, R] = AL_1 \subseteq A$. Therefore, $[w, x_2]x_2 A \subseteq A$ for all $x_2, w \in R$.

Let $N$ be the additive subgroup of $R$ generated by all elements $[w, x_2]x_2$ for $w, x_2 \in R$. Then $NA \subseteq A$. In view of Lemma 7.2, $N \nsubseteq Z(R)$. In view of [6, Theorem, p.98], there exists a proper Lie ideal $L_2$ of $R$ contained in $N$ unless $R \cong \mathrm{M}_2(\mathrm{GF}(2))$.

For the latter case, it follows from Theorem 7.5 that either

$$N = \{0, e_{12} + e_{21}, 1 + e_{12}, 1 + e_{21}\} \text{ or } N = \{0, 1, e_{11} + e_{12} + e_{21}, e_{22} + e_{12} + e_{21}\}.$$

Moreover, $[w, x_2 + 1](x_2 + 1) - [w, x_2]x_2 = [w, x_2] \in N$ for all $w, x_2 \in R$. That is, $[R, R] \subseteq N$. In particular, $e_{12} \in N$. This is a contradiction.

Up to now, we have proved that $L_2 \subseteq N$. Therefore, we get $L_2 A \subseteq A$ and $AL_1 \subseteq A$, where $L_1$ and $L_2$ are proper Lie ideals of $R$. It follows that

$$(L_2 + L_2^2)A \subseteq A \text{ and } A(L_1 + L_1^2) \subseteq A.$$

Set $I_i := R[L_i, L_i]R$ for $i = 1, 2$. In view of Lemma 7.4, $I_1$ and $I_2$ are nonzero. By Lemma 3.11, $I_i \subseteq L_i + L_i^2$ for $i = 1, 2$. Therefore, $I_2 A \subseteq A$ and $AI_1 \subseteq A$. Hence, $0 \neq I_2 A I_1 \subseteq A$, where $I_2 A I_1$ is an ideal of $R$. as desired. □

# 8 Some generalizations

Let $R$ be a prime ring with extended centroid $C$, and let $Q$ denote the Martindale symmetric ring of quotients of $R$. Recall that $Q$ is a prime ring with center $C$. Therefore, $Q$ is an algebra over the field $C$. It is also known that, given $\beta \in C$, there exists a nonzero ideal $I$ of $R$ such that $\beta I \subseteq R$ (see [1, Chapter 2.3] for details). Therefore, if $R$ is a simple ring, then $R = RC$.

**Definition 4.** Given a positive integer $n \geq 2$ and an element $0 \neq \beta \in C$, for $a_1, \ldots, a_n \in R$ we define

$$[a_1, \ldots, a_n]_{n,\beta} := a_1 a_2 \cdots a_n - \beta a_n a_{n-1} \cdots a_1,$$

which is called the $(n, \beta)$-generalized commutator of $a_1, \ldots, a_n$.

Therefore, $[a_1, \ldots, a_n]_n = [a_1, \ldots, a_n]_{n,1}$. We set

$$f_{n,\beta} := [X_1, \ldots, X_n]_{n,\beta} = X_1 X_2 \cdots X_n - \beta X_n X_{n-1} \cdots X_1$$

for $n \geq 2$. It is particularly interesting when $\beta = 1$ or $\beta = -1$. First, we can actually extend Theorem 5.1 as follows.

**Theorem 8.1.** *Let $R$ be a noncommutative prime ring, and let $A$ be a nonzero additive subgroup of $R$, and $0 \neq \beta \in C$, $n \geq 3$ a positive integer. Suppose that $[x_1, \ldots, x_r, a, y_1, \ldots, y_s]_{n,\beta} \in A$ for all $x_i, y_j \in R$ and all $a \in A$, where $r + s = n - 1$. Then $A$ contains a nonzero ideal of $R$.*

We can prove this theorem in a similar way proving Theorem 5.1. Some appropriate modifications are of course necessary. For instance, we need to prove the generalization of Proposition 6.1 (i.e. Theorem 8.1 with $n = 3$ and $r = 1 = s$): Given a nonzero additive subgroup $A$ of a noncommutative prime ring $R$, if $[R, A, R]_{3,\beta} \subseteq A$ then $A$ contains a nonzero ideal of $R$. The necessary modification is the following key identity. Let $a, b \in A$ and $x, y, z \in R$. Then

$$\begin{aligned}
&y(axb - \beta bxa)z \\
&= (yax)bz - \beta zb(yax) + \beta(zby)ax - xa(zby) + (xaz)by - \beta yb(xaz) \\
&= [yax, b, z]_{3,\beta} - [x, a, zby]_{3,\beta} + [xaz, b, y]_{3,\beta} \in A.
\end{aligned}$$

The next step is to prove the following generalization of Lemma 6.4: Let $R$ be a prime ring with extended centroid $C$. If $xay - \beta yax \in Ca$ for all $x, y \in R$, where $0 \neq a \in R$ and $0 \neq \beta \in C$, then $R$ is commutative. We will omit all details in order to keep this paper concise.

Recall that, for $f \in C\{X_1, \ldots, X_n\}$, let $\mathrm{Add}_f(R)$ denote the additive subgroup of $RC$ generated by all elements $f(x_1, \ldots, x_n)$ for $x_i \in R$. We say that $f$ is central-valued on $RC$ if $f(x_1,\dots,x_n)\in C$ for all $x_i\in RC$. Motivated by Corollaries 4.12 and 4.13, it is natural to raise the following.

**Problem 2.** Given a prime ring $R$, characterize polynomials $f\in C\{X_1,\dots,X_t\}$ such that either $\mathrm{Add}_f(R)=R$ or $\mathrm{Add}_f(R)$ contains a nonzero ideal of $R$.

Although we do not understand the general situation, the following theorem provides a partial answer to Problem 2.

**Theorem 8.2.** *Let $R$ be a prime ring, and let*

$$f=\sum_{k=1}^{s} g_k(X_2,\dots,X_n)X_1h_k(X_2,\dots,X_n),$$

*where $g_k,h_k\in C\{X_2,\dots,X_n\}$ for $k=1,\dots,s$. Suppose that $f$ is not central-valued on $RC$ and that $\sum_{k=1}^{s} h_k(X_2,\dots,X_n)g_k(X_2,\dots,X_n)$ is not a PI for $R$. Then $\mathrm{Add}_f(R)$ contains a nonzero ideal of $R$ except when $R\cong \mathrm{M}_2(\mathrm{GF}(2))$ and*

$$\mathrm{Add}_f(R)=\{0,e_{12}+e_{21},1+e_{12},1+e_{21}\}\text{ or }\{0,1,e_{11}+e_{12}+e_{21},e_{22}+e_{12}+e_{21}\}.$$

*Proof.* We assume that the exceptional case is excluded. In view of Theorem 7.5, there exists a nonzero ideal $M$ of $R$ such that $[M,R]\subseteq \mathrm{Add}_f(R)$. Note that, given a nonzero element $\mu\in C$, there exists a nonzero ideal $I$ of $R$ such that $\mu I\subseteq R$. Therefore, we can choose a nonzero ideal $I$ of $R$ contained in $M$ such that $\mathrm{Add}_{g_k}(I)\cup \mathrm{Add}_{h_k}(I)\subseteq M$ for $k=1,\dots,s$.

Let $x_1,\dots,x_n\in I$. Then

$$\begin{aligned}&x_1\textstyle\sum_{k=1}^{s} h_k(x_2,\dots,x_n)g_k(x_2,\dots,x_n)\\ &=\textstyle\sum_{k=1}^{s} g_k(x_2,\dots,x_n)x_1h_k(x_2,\dots,x_n)+\sum_{k=1}^{s}\big[x_1h_k(x_2,\dots,x_n),g_k(x_2,\dots,x_n)\big]\\ &=f(x_1,\dots,x_n)+\textstyle\sum_{k=1}^{s}\big[x_1h_k(x_2,\dots,x_n),g_k(x_2,\dots,x_n)\big].\end{aligned}\tag{9}$$

Note that, in Eq.(9), $\sum_{k=1}^{s}\big[x_1h_k(x_2,\dots,x_n),g_k(x_2,\dots,x_n)\big]\in[M,M]\subseteq \mathrm{Add}_f(R)$. Therefore, $x_1\sum_{k=1}^{s} h_k(x_2,\dots,x_n)g_k(x_2,\dots,x_n)\in \mathrm{Add}_f(R)$ for $x_1,\dots,x_n\in I$. That is,

$$\mathrm{Add}_p(I)\subseteq \mathrm{Add}_f(R),$$

where $q(X_2,\dots,X_n):=\sum_{k=1}^{s} h_k(X_2,\dots,X_n)g_k(X_2,\dots,X_n)$ and

$$p(X_1,\dots,X_n):=X_1q(X_2,\dots,X_n).$$

Since $q(X_2,\dots,X_n)$ is not a PI for $R$ and $R$ and $I$ satisfy the same PIs over $C$ (see [1, Theorem 6.4.1]), it follows that $q(X_2,\dots,X_n)$ is not a PI for $I$. Clearly, $I$ is itself a prime ring with extended centroid $C$. In view of [5, Lemma 5], $\mathrm{Add}_p(I)$ contains a nonzero ideal, say $J$, of the prime ring $I$. Then $IJI$ is an ideal of $R$ and

$$0 \neq IJI \subseteq J \subseteq \mathrm{Add}_p(I) \subseteq \mathrm{Add}_f(R),$$

as desired. □

**Corollary 8.3.** *Let $R$ be a simple ring, and let*

$$f(X_1, \ldots, X_n) = \sum_{k=1}^{s} g_k(X_2, \ldots, X_n) X_1 h_k(X_2, \ldots, X_n),$$

*where $g_k, h_k \in C\{X_2, \ldots, X_n\}$ for $k = 1, \ldots, s$. Suppose that $f$ is not central-valued on $R$ and that $\sum_{k=1}^{s} h_k(X_2, \ldots, X_n) g_k(X_2, \ldots, X_n)$ is not a PI for $R$. Then $\mathrm{Add}_f(R) = R$ except when $R \cong \mathrm{M}_2(\mathrm{GF}(2))$ and*

$$\mathrm{Add}_f(R) = \{0, e_{12} + e_{21}, 1 + e_{12}, 1 + e_{21}\} \text{ or } \{0, 1, e_{11} + e_{12} + e_{21}, e_{22} + e_{12} + e_{21}\}.$$

We next give the following application to Theorem 8.2.

**Theorem 8.4.** *Let $R$ be a noncommutative prime ring and $n \geq 3$ a positive integer, and let $L$ be a Lie ideal of $R$ with $[L, L] \neq 0$. Then $[R_1, \ldots, R_r, L, R_1, \ldots, R_s]_{n,\beta}$ contains a nonzero ideal of $R$, where $R_i = R$ for all $i$ and $r + s = n - 1$.*

*Proof.* By Lemma 3.11, $R[L, L]R \subseteq L + L^2$. Set $I := R[L, L]R$, which is a nonzero ideal of $R$. Let $x_1, x_2 \in L$ and $r \in R$. Then $[x_1x_2, r] = -[rx_1, x_2] - [x_2r, x_1]$. Therefore,

$$[L^2, R] \subseteq [RL, L] + [LR, L] \subseteq L.$$

This implies that $[I, R] \subseteq [L + L^2, R] \subseteq L$. Note that $I$ itself is a noncommutative prime ring with extended centroid $C$ and $\big[[I, I], [I, I]\big] \neq 0$ (see Lemma 7.4). Moreover, every nonzero ideal of $I$ always contains a nonzero ideal of $R$. Therefore, it suffices to prove that

$$\big[I_1, \ldots, I_r, [I, I], I_1, \ldots, I_s\big]_{n,\beta}$$

contains a nonzero ideal of $I$, where $I_i = I$ for all $i$. Set

$$f := \big[X_1, \ldots, X_r, [U, V], Y_1, \ldots, Y_s\big]_{n,\beta},$$

where $X_1, \ldots, X_r, U, V, Y_1, \ldots, Y_s$ are distinct noncommutative indeterminates over $C$. We may assume without loss of generality that $r \geq 1$. Our aim is to prove that $\mathrm{Add}_f(I)$ contains a nonzero ideal of $I$. We now apply Theorem 8.2 to the present case.

First, we claim that $f$ is not central-valued on $IC$. Otherwise, we have

$$\big[x_1, \ldots, x_r, [u,v], y_1, \ldots, y_s\big]_{n,\beta} \in C \tag{10}$$

for all $x_1, \ldots, x_r, u, v, y_1, \ldots, y_s \in I$. Since $I$ and $Q$ satisfy the same GPIs (see [1, Theorem 6.4.1]), replacing $x_i = 1$ for $i > 1$ and $y_j = 1$ for all $j$ in Eq.(10) we get

$$x_1[u,v] - \beta[u,v]x_1 \in C \tag{11}$$

for all $x_1, u, v \in Q$. If $\beta \neq 1$, replacing $x_1 = 1$ we get $(\beta - 1)[Q,Q] \subseteq C$, implying that $Q$ is commutative. This is a contradiction. Hence we get $\beta = 1$. By Eq.(11) we get $\big[Q,[Q,Q]\big] \subseteq C$, implying that $Q$ is commutative. This is a contradiction.

Rewrite $f$ as

$$f = X_1X_2\cdots X_r[U,V]Y_1Y_2\cdots Y_s - \beta Y_sY_{s-1}\cdots Y_1[U,V]X_rX_{r-1}\cdots X_1.$$

We claim that

$$g := X_2\cdots X_r[U,V]Y_1Y_2\cdots Y_s - \beta Y_sY_{s-1}\cdots Y_1[U,V]X_rX_{r-1}\cdots X_2$$

is not a PI for $I$. Otherwise, we have

$$x_2\cdots x_r[u,v]y_1y_2\cdots y_s - \beta y_sy_{s-1}\cdots y_1[u,v]x_rx_{r-1}\cdots x_2 = 0 \tag{12}$$

for all $x_2, \ldots, x_r, u, v, y_1, \ldots, y_s \in Q$.

Replacing all $x_i, y_j$ by 1 in Eq.(12), we get $(\beta - 1)[Q,Q] = 0$, implying that $\beta = 1$. If $r \geq 2$, then replacing $x_i = 1 = y_j$ for $i > 2$ and $j \geq 1$ we get $[x_2,[u,v]] = 0$ for all $x_2, u, v \in Q$, implying that $Q$ is commutative, a contradiction. Suppose next that $r = 1$. Then $s \geq 1$. Replacing $y_j = 1$ for $j > 1$, we get $[y_1,[u,v]] = 0$ for all $y_2, u, v \in Q$. As above, we get a contradiction.

In view of Theorem 8.2, $\mathrm{Add}_f(I)$ contains a nonzero ideal of $I$ except when $I \cong \mathrm{M}_2(\mathrm{GF}(2))$ and

$$\mathrm{Add}_f(I) = \{0, e_{12} + e_{21}, 1 + e_{12}, 1 + e_{21}\} \text{ or } \{0, 1, e_{11} + e_{12} + e_{21}, e_{22} + e_{12} + e_{21}\}.$$

We are done if $\mathrm{Add}_f(I)$ contains a nonzero ideal of $I$. Hence we assume that the latter case occurs.

**Case 1**: $r \geq 2$ or $s \geq 2$. We may assume that $r \geq 2$. The case that $s \geq 2$ has the same argument. Replacing $X_1$, $X_2$ by $e_{11}$, $e_{12}$, respectively, and $X_i, Y_j$ by 1 for $i > 2$ and all $j$, $U$ by $e_{22}$, and $V$ by $e_{21}$, we get $e_{11} \in \mathrm{Add}_f(I)$, a contradiction.

**Case 2**: $r = 1 = s$. In this case, we have

$$f(X_1, U, V, Y_1) = X_1[U,V]Y_1 - \beta Y_1[U,V]X_1.$$

Then $f(e_{11}, e_{11}, e_{12}, e_{22}) = e_{12} \in \mathrm{Add}_f(I)$, a contradiction. □

Recall that $f_{n,\beta} := [X_1, \ldots, X_n]_{n,\beta}$ for $n \geq 2$. As an immediate consequence of Theorem 8.1, we have the following corollary.

**Corollary 8.5.** *Let $R$ be a noncommutative prime ring and $n \geq 3$. Then* $\mathrm{Add}_{f_{n,\beta}}(R)$ *contains a nonzero ideal of $R$. In addition, if $R$ is a simple ring, then* $\mathrm{Add}_{f_{n,\beta}}(R) = R$.

**Corollary 8.6.** *Let $R$ be a noncommutative prime ring and $n \geq 4$ a positive integer, and let $L$ be a nonzero Lie ideal of $R$. Then $[R_1, \ldots, R_r, L, R_1, \ldots, R_s]_{n,\beta}$ contains a nonzero ideal of $R$, where $R_i = R$ for all $i$ and $r + s = n - 1$, except when* $\mathrm{char}\, R = 2$ *and* $\dim_C RC = 4$.

*Proof.* We exclude the exceptional case. By Theorem 8.4, we are done if $[L, L] \neq 0$. Suppose next that $[L, L] = 0$. In view of [18, Lemma 7], $L \subseteq Z(R)$. Therefore,

$$[R_1, \ldots, R_r, L, R_1, \ldots, R_s]_{n,\beta} = L[R, \ldots, R]_{n-1,\beta}.$$

Note that $n - 1 \geq 3$. It follows from Corollary 8.5 that $[R, \ldots, R]_{n-1,\beta}$ contains a nonzero ideal of $R$, so does $L[R, \ldots, R]_{n-1,\beta}$, as desired. □

**Acknowledgments**. The work of P. V. Danchev was supported in part by the Bulgarian National Science Fund under Grant KP-06 No. 32/1 of December 07, 2019, and that of T.-K. Lee was supported in part by the Ministry of Science and Technology of Taiwan (MOST 109-2115-M-002-014).